\documentclass[11pt, a4paper]{article}
\usepackage{amsfonts, amsthm, amsmath}
\usepackage{graphicx, geometry, microtype}
\usepackage{amssymb}
\usepackage{epstopdf}

\allowdisplaybreaks[4]

\newtheorem{thm}{Theorem}

\newtheorem{Def}{Definition}
\newtheorem{cor}{Corollary}

\newcommand{\reals}{\mathbb{R}}

\newcommand{\CC}{\mathbb{C}}
\newcommand{\disk}{\mathbb{D}}
\newcommand{\PP}{\mathbb{F}^2}

\newcommand{\Euclideanplane}{\mathbb{E}^2}
\newcommand{\hyperbolicspace}{\mathbb{H}^3}
\newcommand{\hyperbolicplane}{\mathbb{H}^2}
\newcommand{\cellcomplex}{\mathcal{C}}

\newcommand{\Hplane}{\mathbb{H}_0}

\newcommand{\Triang}{\mathcal{T}}

\newcommand{\polytope}{\mathcal{P}_{S,\cellcomplex}}
\newcommand{\naturals}{\mathbb{N}}
\newcommand{\ADelta}{\mathcal{A}_{\Delta}}

\newcommand{\Vol}{\mathbf{V}}

\newcommand{\ER}{\mathcal{ER}}
\newcommand{\TED}{\mathcal{TE}_{\Delta}}
\newcommand{\ERD}{\mathcal{ER}_{\Delta}}
\newcommand{\TE}{\mathcal{TE}}
\newcommand{\UDelta}{\mathbf{U}_{\Delta}}
\newcommand{\Uglobal}{\mathbf{U}}

\date{}

\begin{document}

\title{Discrete uniformization of
finite branched covers over the Riemann sphere via hyper-ideal
circle patterns}

\author{Alexander Bobenko \and Nikolay Dimitrov \and Stefan Sechelmann}

\maketitle

\begin{abstract}
With the help of hyper-ideal circle pattern theory, we have
developed a discrete version of the classical uniformization
theorems for surfaces represented as finite branched covers over
the Riemann sphere as well as compact polyhedral surfaces with
non-positive curvature. We show that in the case of such surfaces
discrete uniformization via hyper-ideal circle patterns always
exists and is unique. We also propose a numerical algorithm,
utilizing convex optimization, that constructs the desired
discrete uniformization.
\end{abstract}

\section{Introduction}

In the current paper, we have constructed discrete conformal maps
by applying the theory of \emph{hyper-ideal circle patterns}.
Informally speaking, a circle pattern is a discrete collection of
overlapping circles together with intersection angles between
adjacent circles, following the combinatorics of a given polygonal
cell complex on a polyhedral surface
\cite{ThuNotes,Riv94,BobSpr1}. In the case of \emph{hyper-ideal
circle patterns}, the collection is divided into two subsets --
one assigned to the vertices of the complex and one assigned to
the faces, so that the face circles are orthogonal to their vertex
neighbors \cite{Sch1,S,ND}. We have chosen this theory because of
its existence and uniqueness theorems which hold true for fixed
underlying combinatorial cell complexes, contrary to some other
discrete conformal theories. The conceptual geometric rationale
behind our choice is based on the combination of two equivalent
characterizing properties of smooth conformal maps. First, at
every point, a smooth conformal map stretches or shrinks the
surface equally in all directions. Hence, at each point
infinitesimal circles are mapped to infinitesimal circles. Second,
smooth conformal maps preserve angles. Since circles and angles
are the main ingredients of circle patterns, one can discretize
the smooth theory by promoting a finite number of infinitesimal
circles into circles on the surface, while keeping track of the
intersection angles between pairs of neighboring ones.

The idea of using discrete collections of circles as a discrete
model of conformal transformations comes from Thurston's
conjecture \cite{ThuAdr} that the Riemann mapping of a simply
connected planar region can be approximated by a sequence of finer
and finer regular, hexagonal circle packings \cite{Steph} (a
collection of touching circles with the combinatorics of the
hexagonal grid on the plane). The conjecture was proved by Rodin
and Sullivan \cite{R_and_S}. After that, circle packings with more
general combinatorics were used to define various discrete
analogues of holomorphic maps on the plane \cite{Steph} as well as
discrete analogue of the uniformization theorem
\cite{Beardon_Steph}. Circle packings however are determined by
combinatorics. In order to allow for more freedom, one can turn
towards standard (Delaunay) circle patterns which provide the
opportunity to incorporate the intrinsic geometry of the surfaces
one works with, and more precisely, the intersection angles
between adjacent circles. For example, some holomorphic
transformations in the complex plane have been discretized using
orthogonal circle patterns with the combinatorics of the square
grid \cite{Schramm,BobSur}, exploring links with integrable
systems. Nevertheless, when handling polyhedral surfaces with cone
singularities, standard circle patterns do not posses the
necessary degrees of freedom that allow manipulation of the cone
angles of the underlying geometry. This obstacle is evident in
\cite{KSS}, where an attempt is made to construct discrete
conformal maps via Delaunay circle patterns. The problem with such
approach is that the definition of conformal equivalence does not
actually form an equivalence relation but is more of an
approximation. In contrast, \emph{Hyper-ideal circle patterns}
\cite{Sch1,S} possess the necessary flexibility and, as we have
shown in the current article, can be used to construct well
structured discrete uniformizing conformal maps for surfaces with
non-positive cone singularities as well as surfaces given as
ramified covers over the sphere.

The classical uniformization theorem is a central result in the
theory of complex analytic functions, conformal maps and Riemann
surface theory. The uniformization theorem states that every
simply connected Riemann surface is conformally isomorphic to
either the Riemann sphere $\hat{\CC}$, the plane $\CC$ or the unit
disc $\disk$, uniquely up to conformal automorphisms. In the case
of a closed Riemann surface of genus one or greater, the theorem
implies that, up to conformal automorphisms, its universal
covering map can be seen as a unique conformal transformation from
either the plane $\CC$ (in the case of genus one) or from the unit
disc $\disk$ (in the case of genus two or higher) with a
fundamental covering group represented by a subgroup of the
conformal automorphism group of $\CC$ or $\disk$ respectively.
Since $\CC$ has a natural Euclidean metric and $\disk$ is a
conformal model of the hyperbolic plane $\hyperbolicplane$ with a
natural metric of constant Gaussian curvature $-1$, the
uniformization theorem can be reinterpreted as the fact that each
closed Riemann surface has a natural homogeneous metric of
constant Gaussian curvature.

Consequently, when discretizing the uniformization theorem, one
could attempt to obtain the same properties for the discrete
uniformizing map as the ones discussed above. Also, in general,
the goal should be to start from a polyhedral metric (or some
other general structure that determines a conformal structure) on
the surface and obtain the discretely conformally equivalent
smooth metric of constant curvature. For instance in
\cite{Beardon_Steph}, a circle packing approximation of the
universal covering map of a Riemann surface is presented, but the
question is more about inducing a metric with a circle packing on
a topological triangulated surface with given combinatorics by
laying it out on the universal cover. In \cite{BPS,GLSW1,GLSW2} a
flexible and constructive discrete uniformization is obtained
based on the method of \emph{discrete conformal factors} for
triangulated surfaces \cite{Luo,SSP}. It is based on the discrete
analog of conformally equivalent metrics. In the case of
polyhedral surfaces with fixed combinatorics, there is a
uniqueness theorem but there is no existence theorem. In other
words, it is not always guaranteed that the method will produce a
discretely conformally equivalent image if the combinatorics is
preserved. However, if one is allowed to perform edge-flips and
thus alter the combinatorial nature of the triangular mesh, a
solution is guaranteed \cite{GLSW1,GLSW2}.

Another interesting approach to discrete conformal mappings is the
method of \emph{inverse distance circle packings}. It has been
introduced in \cite{BowSteph} and tested as a possible
discretization of conformal mappings in \cite{BowHurd}. Inverse
distance packings are somewhat similar to hyper-ideal circle
patterns, except that instead of intersection angles between
adjacent face circles, inverse distances between vertex-circles
are prescribed. In \cite{LuoIII}, Luo has proved a uniqueness
theorem (i.e. a rigidity theorem) for inverse distance packings in
the Euclidean and hyperbolic cases. There is no uniqueness in the
spherical case since a counter-example has been constructed in
\cite{JimSchl}. No existence theorem for inversive distance
packings is known.

Our discrete uniformization method describes a discrete conformal
map in terms of hyper-ideal circle patterns (see Definition
\ref{Def_discrete_conformal_equivalence}). It features a
well-defined conformal equivalence relation and has the desired
properties. In contrast to the methods of conformal factors and
inverse distance packings, it possesses the existence and
uniqueness theorems for fixed polygonal combinatorics.
Furthermore, it is constructive and it has been algorithmically
implemented. In the last section we have included some computer
generated examples of discrete uniformization via hyper-ideal
circle patterns. However, our approach is restricted to surfaces
with polyhedral metrics of non-positive curvature, i.e. the cone
singularities have angles greater than $2\pi$,  as well as to
finite ramified covers over the sphere. All algebraic curves fall
in the latter category.






\section{Definitions and notations} \label{Sec_def_and_notations}

We start with some terminology and notations. Assumed that $S$ is
an orientable compact topological surface with no boundary. Denote
by $d$ a metric of constant negative or zero Gaussian curvature on
$S$ with finitely many cone singularities $\text{sing}(d)$.
Depending on whether the curvature away from singularities is
negative or zero, the metric will be referred to as a
\emph{hyperbolic cone-metric} or a \emph{Euclidean cone-metric}
respectively. We will use $\mathbb{F}^2$ as a common notation for
both the Euclidean plane $\Euclideanplane$ and the hyperbolic
plane $\hyperbolicplane$. In this paper, we mostly use the
Poincar\'e disc model and the upper half-plane model of
$\hyperbolicplane$. Both of them are conformal and they have the
advantage that circles in hyperbolic geometry appear in both
models as circles in the underlying Euclidean geometry (for
details, see \cite{ThuNotes,ThuBook,BenPetr}). The only
particularity is that, in general, the hyperbolic centers of the
circles do not necessarily coincide with their Euclidean ones.
Throughout this article, we will also use the notation $\hat{\CC}
= \CC \cup \{\infty\}$ to denote the Riemann sphere, which could
also be thought of as the projective complex line
$\CC\mathbb{P}^1$. The global conformal automorphisms of
$\hat{\CC} \cong \CC\mathbb{P}^1$ are the M\"obius transformations
(i.e. the linear fractional transformations) which form the group
$\mathbb{P}SL(2,\CC)$. From now on, $V$ will be a finite set of
points on $S$ containing the cone singularities of $d$. Thus $V$
splits into two disjoint subsets $V_1 = \text{sing}(d)$ and $V_0 =
V \setminus V_1$. By $\cellcomplex=(V,E,F)$ we will denote a
topological cell complex on $S$, where $V$ are its vertices,
partitioned into two disjoint subsets $V_0$ and $V_1$, $E$ are the
edges and $F$ are the faces of $\cellcomplex$ (see for example
Figure \ref{Fig2} below). All three sets are finite. Furthermore,
without loss of generality, we will always assume that the cell
complexes we work with, their dual complexes, and the various
subdivisions of the former and the latter, are nicely embedded in
the surface. Moreover, they will be assumed to be always
\emph{strongly regular} \cite{Z,G,S1} which means that (i) closed
cells (of any dimension) are attached without identifications on
their boundaries and (ii) the intersection of any pair of closed
cells is either a closed cell or empty.


\begin{Def} \label{Def_geodesic_complex}
A geodesic cell complex on $(S,d)$ is a cell complex
$\cellcomplex_d=(V,E_d,F_d)$ whose edges, with endpoints removed,
are open geodesic arcs embedded in $S \setminus V$. 
\end{Def}

In other words, we can think of a geodesic cell-complex
$\cellcomplex_d$ on a geometric surface $(S,d)$ as a two
dimensional manifold, obtained by gluing together geodesic
polygons along their edges. The edges that are being identified
should have the same length and the identification should be an
isometry. Observe the difference between a geodesic cell-complex
$\cellcomplex_d$ and its topological (combinatorial) counterpart
$\cellcomplex$. While $\cellcomplex_d$ is made of geodesic
polygons and thus provides the underlying surface $S$ with a
cone-metric $d$, the cell-complex $\cellcomplex$ is just a purely
topological (and hence combinatorial) object. In many cases in
this article $\cellcomplex$ will be obtained from $\cellcomplex_d$
by forgetting about the geometry of $\cellcomplex_d$ and focusing
entirely on its combinatorics and topology.

Assume three disjoint circles $c_i, c_j$ and $c_k$ with centers
$i, j$ and $k$ respectively, lie in the geometric plane $\PP$.

\begin{figure}
\centering
\includegraphics[width=8.7cm]{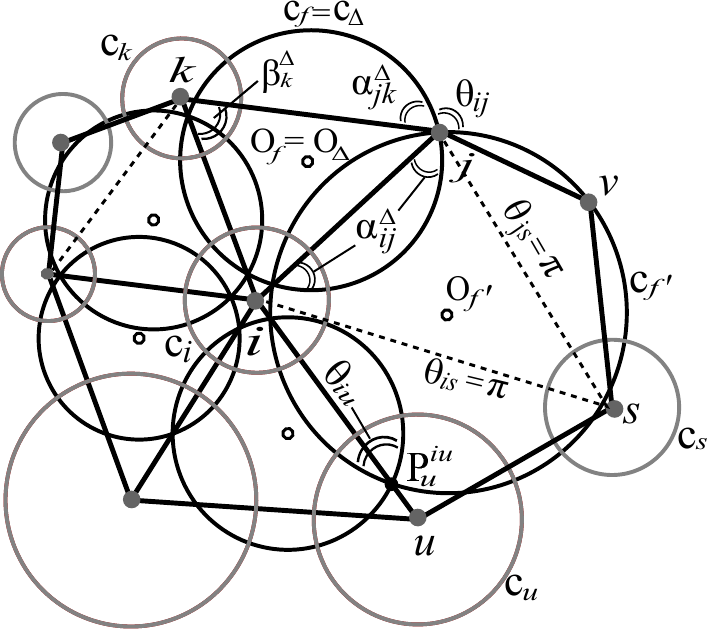}
\caption{A hyper-ideal circle pattern of decorated polygons
together with labels.} \label{Fig1}
\end{figure}

Then, there exists a unique forth circle $c_{\Delta}$ orthogonal
to $c_i, c_j$ and $c_k$. Furthermore, let $\Delta=ijk$ be the
geodesic triangle spanned by the centers $i, j$ and $k$. Then
$\Delta$, together with the circles $c_i, c_j, c_k$ and
$c_{\Delta}$, is called a \emph{decorated triangle} (see Figure
\ref{Fig1}). The circles $c_i, c_j$ and $c_k$ are called the
\emph{vertex circles} of $\Delta$, while $c_{\Delta}$ is called
the \emph{face circle} of $\Delta$. We emphasize that it is
allowed for one, two or all three vertex circles to have zero
radii, i.e. to be points. Even in this more general set up,
everything said above still applies. Whenever the vertex circles
have zero radius, we will call them \emph{degenerate vertex
circles}.

Now, assume two non-overlapping decorated triangles, like
$\Delta_1=jis$ and $\Delta_2=uis$ from Figure \ref{Fig1}, share a
common edge $is$. As usual, denote by $c_i, c_j, c_s$ and $c_u$
the vertex circles (some of which may be degenerate), and by
$c_{\Delta_1}$ and $c_{\Delta_2}$ the corresponding face circles
of the triangles. Although generically the two face circles
$c_{\Delta_1}$ and $c_{\Delta_2}$ are different, sometimes it may
happen that they coincide, i.e. $c_{\Delta_1}=c_{\Delta_2}=c_q$.
In that case all four vertex circles $c_i, c_j, c_s$ and $c_u$ are
orthogonal to $c_q$. Thus, we can erase the edge $is$ and obtain a
decorated geodesic quadrilateral $q=ijsu$ with vertex circles
$c_i, c_j, c_s$ and $c_u$, and a face circle $c_q$ orthogonal to
the vertex ones. Observe, that in this case the quadrilateral is
convex. If we continue this way, we can obtain decorated polygons
with arbitrary number of edges, like for instance the decorated
pentagon $f'=ijvsu$ from Figure \ref{Fig1}.

\begin{Def}\label{Def_decorated_polygon}
A decorated polygon is a convex geodesic polygon $p$ in $\PP$,
with vertices $i_1, i_2, ..., i_n$, together with:

\smallskip

\noindent $\bullet$ a set of disjoint circles $c_{i_1}, c_{i_2},
..., c_{i_n}$ such that each $c_{i_s}$ is centered at the vertex
$i_s$ for $s=1,..,n$. Some or all of the circles are allowed to be
degenerate, i.e. circles of radius zero;

\smallskip

\noindent $\bullet$ another circle $c_p$ orthogonal to
$c_{i_1},...,c_{i_n}$.

\smallskip

\noindent The circles $c_{i_1},...,c_{i_n}$ are called
\emph{vertex circles} and the additional orthogonal circle $c_p$
is called the \emph{face circle} of the decorated polygon $p$. 
\end{Def}

Two faces of a cell complex that share a common edge will be
called \emph{adjacent to each other}. Furthermore, assume two
decorated polygons $p_1$ and $p_2$ share a common geodesic edge
$ij$. Then, the decorated polygons $p_1$ and $p_2$ are called
\emph{compatibly adjacent} whenever the vertex circles $c_i^1,
c_j^1$ of $p_1$ and $c_i^2, c_j^2$ of $p_2$ coincide respectively,
that is $c_i^1\equiv c_i^2$ and $c_j^1\equiv c_j^2$. Furthermore,
whenever two decorated polygons are compatibly adjacent, we will
say that their face circles are \emph{adjacent to each other}. A
situation like that is depicted on Figure \ref{Fig1} for the edge
$ij$ and the two faces $f\equiv\Delta$ and $f'$ with face circles
$c_f$ and $c_{f'}$.

\begin{Def} \label{Def_local_Del_property}
Let $p_1$ and $p_2$ be two compatibly adjacent decorated polygons
in $\PP$. Let $ij$ be their common geodesic edge. Furthermore, let
$c_{p_1}$ and $c_{p_2}$ be the face circles of $p_1$ and $p_2$
respectively.

\smallskip

\noindent $\bullet$ We say that the edge $ij$ satisfies the
\emph{local Delaunay property} whenever each vertex circle of the
decorated polygon $p_2$ is either (i) disjoint from the interior
of the face circle $c_{p_1}$ of $p_1$, or (ii) if it is not, the
intersection angle between the vertex circle in question and the
face ci rcle $c_{p_1}$, measured between the circular arcs that
bound the region of common intersection of their discs, is less
than $\pi /2$. See for instance edge $ij$ on Figure \ref{Fig1}.

\smallskip

\noindent $\bullet$ For the edge $ij$, which satisfies the local
Delaunay property, $\theta_{ij} \in (0,\pi)$ denotes the
intersection angle between the two adjacent face circles $c_{p_1}$
and $c_{p_2}$, measured between the circular arcs that bound the
region of common intersection of their discs. (See for example
angles $\theta_{ij}$ and $\theta_{iu}$  
from Figure \ref{Fig1}.)
\end{Def}

It is not difficult to see that the definition of local Delaunay
property is symmetric in the sense that if the condition of
Definition \ref{Def_local_Del_property} holds for the face circle
$c_{p_1}$ and the vertex circles of $p_2$, then it also holds for
the face circle $c_{p_2}$ and the vertex circles of $p_1$.

\begin{Def} \label{Def_hyperideal_circle_pattern}
A \emph{hyper-ideal circle pattern} on a given surface $S$ (Figure
\ref{Fig1}) is a hyperbolic or Euclidean cone metric $d$ on $S$
together with a geodesic cell complex $\cellcomplex_d=(V,E_d,F_d)$
whose faces are decorated geodesic polygons such that any two
adjacent faces are compatibly adjacent and each geodesic edge of
$\cellcomplex_d$ has the local Delaunay property. Whenever $d$ is
flat on $S\setminus V$, we call the circle pattern Euclidean, and
whenever $d$ is hyperbolic on $S\setminus V$, we call the pattern
hyperbolic. Whenever all vertex radii are zero, the pattern will
be called a \emph{Delaunay circle pattern} (see left-hand side of
Figure \ref{Fig2}).
\end{Def}

Intuitively speaking, a hyper-ideal circle pattern on a surface
$S$ is a surface homeomorphic to $S$, obtained by gluing together
decorated geodesic polygons along pairs of corresponding edges.
The edges that are being identified should have the same length,
the identification should be an isometry and the vertices that get
identified should have vertex-circles with same radii.


Observe that a hyper-ideal circle pattern on $S$ consists of (i) a
cone-metric $d$ on $S$, (ii) a finite set of vertices $V \supseteq
\text{sing}(d)$, (iii) an assignment of vertex radii $r$ on $V$,
and (iv) a geodesic cell complex $\cellcomplex_d$ together with
(v) a collection of vertex circles and (vi) a collection of face
circles. However, the geometric data $(S, d, V, r)$ is enough to
further identify uniquely the geodesic cell complex
$\cellcomplex_d$ and the collections of vertex and face circles.
This is done via the weighted Delaunay cell decomposition
construction \cite{EB}. More precisely, given (i) a geometric
surface $(S,d)$, (ii) a finite set of points $V \supset
\text{sing}(d)$ on $S$ and (iii) an assignment of disjoint vertex
circle radii $r : V \to [0,\infty)$, one can uniquely generate
(obtain) the corresponding $r-$weighted Delaunay cell complex
$\cellcomplex_d$, where each edge satisfies the local Delaunay
property. In the process, the families of vertex and face circles
naturally appear as part of the construction \cite{EB,Sch1,S}.
Thus, hyper-ideal circle patterns are the same as a certain type
of weighted Delaunay cell decompositions. Generically, a weighted
Delaunay cell complex is in fact a triangulation, also known as a
weighted Delaunay triangulation \cite{EB}.

\paragraph{Discussion on planar weighted Delaunay tessellations.} As a
brief illustration, we discuss the latter in the case when the
geometric surface $S,d$ is the plane $\PP$. So we are given a
finite set of points $V$ in $\PP$ together with weights $r : V \to
[0,\infty)$. This is equivalent to actually having a finite set of
closed circular (vertex) disks in $\PP$ with centers $V$ and radii
$r$. Assume they are disjoint. As already discussed in the
paragraphs preceding Definition \ref{Def_decorated_polygon}, for
every triple of such disks there exists a circle that intersects
the boundaries of the disks orthogonally. The weighted Delaunay
triangulation, induced by the vertex disks, consists of those
geodesic triangles whose vertices are the centers of triples of
disks for which the orthogonal circle of this triple intersects no
other disk more than orthogonally (see Figure \ref{Fig1} when
$\PP=\Euclideanplane$ and the right-hand side of Figure \ref{Fig2}
when $\PP=\hyperbolicplane$). In some cases, an orthogonal circle
can intersect orthogonally more than three circles, thus forming a
polygon rather than a triangle, resulting in a more general cell
complex than a triangulation, usually referred to as a
\emph{weighted Delaunay tessellation} and a \emph{weighted
Delaunay cell complex} (an example depicted on Figures \ref{Fig1}
and \ref{Fig2}). Consequently, the vertex circles and the
orthogonal (face) circles form a hyper-ideal circle pattern on
$\PP$ with convex geodesic boundary. The boundary geodesic edges
are technically also circles which in some cases, such as the case
of $\PP=\Euclideanplane$ for instance, pass through the point of
infinity and in the case of $\PP=\hyperbolicplane$ are orthogonal
to the ideal boundary of the hyperbolic plane, which can be
thought of as a vertex circle centered at the point of infinity.
Alternatively, one can obtain the $r-$weighted Delaunay cell
complex as the geodesic dual to the $r$-weighted Voronoi diagram
\cite{EB}. A Voronoi cell in the Euclidean case is defined as
$W_{r}(i)=\big\{ x \in \Euclideanplane \, | \,
d_{\Euclideanplane}(x,i)^2-r_i^2 \leq
d_{\Euclideanplane}(x,j)^2-r_j^2 \,\, \text{for all} \,\, j \in V
\, \big\}$, and a Voronoi cell in the hyperbolic case is defined
as $W_{r}(i)=\big\{ x \in \hyperbolicplane \, | \,
\cosh{(r_j)}\cosh{d_{\hyperbolicplane}(x,i)} \leq
\cosh{(r_i)}\cosh{d_{\hyperbolicplane}(x,j)} \,\, \text{for all}
\,\, j \in V \, \big\}$. The weighted Voronoi complex and the
weighted Delaunay complex are dual to each other and are geodesic
and embedded. The vertices of the Voronoi complex are the centers
of the Delaunay (face) circles of the pattern. For a rough
illustration of a complex and its dual, see Figure \ref{Fig3}.
Both the weighted-Delaunay construction and the Voronoi
construction can be performed on the surface $S,d$ in an analogous
manner.

\medskip

Intuitively speaking, given a topological surfaces $S$ with a
finite number of points on it $V$, one can introduce a discrete
conformal structure on $(S,V)$ by assigning (i) either a
hyperbolic or Euclidean cone metric $d$ such that $\text{sing}(d)
\subseteq V$ together with (ii) an appropriate vertex radii
assignment $r : V \to [0, \infty)$. In short, $(S,V,d,r)$ could be
regarded as a surface with a discrete conformal structure, i.e. a
discrete Riemann surface. Alternatively, instead of a cone-metric
on $S$, one could also have a projective (i.e. conformal)
$\CC\mathbb{P}^1 \cong \hat{\CC}$ structure with cone
singularities, where the latter are assumed to be among the points
from $V \subset S$. What we mean is that away from the points from
$V_1 \subset V$, the surface has an atlas with transition
functions given exclusively by M\"obius transformations from the
conformal group $\mathbb{P}SL(2,\CC)$. Since the elements of
$\mathbb{P}SL(2,\CC)$ send circles to circles and preserve angles
on $\hat{\CC}$, both the notions of circles and angles (but not
centers and radii of circles) are well defined on the surface
itself. Consequently, hyper-ideal circle patterns on a surface
with $\hat{\CC}$ cone-structure make perfect sense and can be
regarded as discrete conformal structures. Notice that the metric
definitions above (Definition \ref{Def_hyperideal_circle_pattern})
are special cases of such cone $\hat{\CC}$ structures because both
the Euclidean and the hyperbolic isometry groups are subgorups of
$\mathbb{P}SL(2,\CC)$.

\paragraph{Discussion on hyper-ideal circle patterns on $\hat{\CC}$.} In
this article, somewhat implicitly, we will encounter Delaunay
circle patterns (but not hyper-ideal ones) on surfaces with cone
$\CC\mathbb{P}^1$ structures (Section
\ref{Sec_uniformization_hyper_elliptic_surfaces}). The only
explicit Delaunay and hyper-ideal circle patterns with $\hat{\CC}$
structure we actually encounter in this paper, are the ones on
$\hat{\CC}$ itself. They can be acted upon by
$\mathbb{P}SL(2,\CC)$, which preserves the combinatorics and the
intersection angles, but does not preserve the notion of circle
centers and radii. Hence, patterns on $\hat{\CC}$ are defined up
to $\mathbb{P}SL(2,\CC)$ M\"obius transformations and for that
reason the circles from such patterns do not have naturally
defined centers and the patterns themselves do not have
canonically defined geometric cell complexes on $\hat{\CC}$
associated to them, although they do have canonically defined
underlying combinatorial complexes. However, Delaunay
\cite{BobSpr1,Riv94,ThuNotes} and hyper-ideal circle patterns
\cite{BaoBon,Rou,Schl} naturally correspond to convex ideal and
hyper-ideal hyperbolic polyhedra in $\hyperbolicspace$
respectively (e.g. see
 as well as Section
\ref{Sec_Spaces_of_triangles_and_patterns}). These polyhedra have
natural geodesic cell decompositions, determined by their
vertices, faces and edges, and represent the combinatorial cell
complexes arising from the circle patterns. Alternatively, one can
directly use Definition \ref{Def_hyperideal_circle_pattern}
without resorting to hyperbolic polyhedra. For example, this can
be done by performing a stereographic projection of the
hyper-ideal circle pattern on $\hat{\CC}$ to the plane and thus
obtaining a hyper-ideal circle pattern with polygonal boundary
either in the Euclidean or in the hyperbolic plane, just as
described in the Discussion on planar weighted Delaunay
tessellations two paragraphs above (see also Figures \ref{Fig1}
and \ref{Fig2}). Whenever one performs a stereographic projection
from a vertex circle shrunk to a point, one obtains a Euclidean
hyper-ideal circle pattern with convex polygonal boundary.
Alternatively, one can also pick a point from the interior of a
vertex circle (the interior disjoint from all other vertex
circles) and then the stereographic projection with respect to
that point leads to a hyper-ideal circle pattern with convex
polygonal geodesic boundary in $\hyperbolicplane$, where $\partial
\hyperbolicplane$ is the image of the vertex circle whose interior
contains the point. As a result, one obtains a geodesic
representation of the cell complex defining the combinatorics of
the pattern in the convex geodesic polygon, as depicted on Figures
\ref{Fig1} and \ref{Fig2}. To complete the cell decomposition to
the sphere, one can add (i) in the Euclidean case straight rays
from the vertices of the boundary going away to the point at
infinity, for example following the angle bisectors at the
vertices; (ii) in the hyperbolic case one can draw the straight
(Euclidean) rays connecting the vertices of the hyperbolic
boundary polygon to the point at infinity, each ray starting from
a vertex and passing through its inverse image with respect to
$\partial \hyperbolicplane$ (when the polygon is in general
position).

\paragraph{Remark:} Just like in the Euclidean and the hyperbolic
cases, one could define hyper-ideal circle patterns with
underlying cone metrics of constant positive curvature (spherical
cone metrics). For instance, one could fix such a metric on
$\hat{\CC}$ and use the metric Definition
\ref{Def_hyperideal_circle_pattern} to define hyper-ideal circle
patterns with a spherical metric on $\hat{\CC}$ and even on $S$
(with cone singularities). Consequently, one can obtain a
spherical geodesic Delaunay cell complex and its Voronoi dual
representing the combinatorics of the pattern. However, it seems
more natural to use the conformal structure on the Riemann sphere,
rather than just one fixed metric. For that reason, we are going
to skip the spherical case in this article.

Given a hyper-ideal circle pattern, one can extract from it the
finite combinatorial data $(\cellcomplex, \theta, \Theta)$, where
(i) $\cellcomplex$ is the geodesic cell-complex $\cellcomplex_d$
on $(S,d)$, regarded as a purely combinatorial object (i.e. we
forget all geometric information), (ii) $\theta \, : \, E \, \to
\, (0, \pi)$ is the assignment of intersection angles of all pairs
of adjacent face circles of the pattern and (iii) $\Theta \, : \,
V \, \to \, (0,+\infty)$ are the cone angles around the points
from $V$. In this case we will say that the hyper-ideal circle
pattern \emph{realizes the (combinatorial) angle data}
$(\cellcomplex, \theta, \Theta)$. As pointed out in the definition
of the set $V$, the cone angle $\Theta_k = 2 \pi$ whenever $k \in
V_0$ and $\Theta_k \neq 2 \pi$ when $k \in V_1 = \text{sing}(d)$.

\begin{Def}\label{Def_discrete_conformal_equivalence}
Two hyper-ideal circle patterns on $S$ are considered
\emph{discretely conformally equivalent} whenever their underlying
(geodesic) cell-complexes are combinatorially isomorphic and the
corresponding intersection angles between pairs of adjacent face
circles are equal. Consequently, a \emph{discrete conformal map}
between two conformally equivalent hyper-ideal circle patterns is
the pairing of face and vertex circles from the first pattern with
the corresponding face and vertex circles from the second pattern.
\end{Def}

\section{Discrete uniformization of polyhedral surfaces with non-positive
curvature.} \label{Sec_uniformization_discrete_negative_surfaces}

Consider the following input data on the surface $S$:

\smallskip

\noindent (a) either Euclidean or hyperbolic cone-metric $d$ on
$S$;

\smallskip
\noindent (b) a finite set of point $V \subset S$ such that $V =
V_0 \cup V_1$, where $V_1 = \text{sing}(d)$ and $V_0 \cap V_1 =
\varnothing$;

\smallskip
\noindent (c) the cone angle $\Theta_k > 2\pi$ when $k \in V_1$
and $\Theta_k = 2\pi$ when $k \in V_0$.

\smallskip
\noindent To put it shortly, the surface $S$ is provided with a
cone metric of non-positive Gaussian curvature. We will use the
notation $(S,d,V)$ to denote the aforementioned input data. One
could interpret it as the geometric data $(S,d,V,0)$ that gives
rise to the classical Delaunay circle pattern, where all vertex
radii are equal to zero.

Next, from the geometry of $(S,d,V)$, form the combinatorial angle
data $(\cellcomplex, \theta, 2 \pi)$, where

\smallskip
\noindent 1. $\cellcomplex = (V,E,F)$ is the unique Delaunay cell
complex of $(S,d)$ with respect to $V$, regarded as a purely
topological (combinatorial) complex.

\smallskip
\noindent 2. $\theta \, : \, E \to \, (0,\pi)$  are the angles
between the pairs of adjacent Delaunay circles of the Delaunay
circle pattern.


\smallskip
\noindent 3. $\Theta_k = 2\pi$ for all $k \in V$. We use the
notation $2\pi$ to denote this constant angle assignment.

\begin{figure}
\centering
\includegraphics[width=14cm]{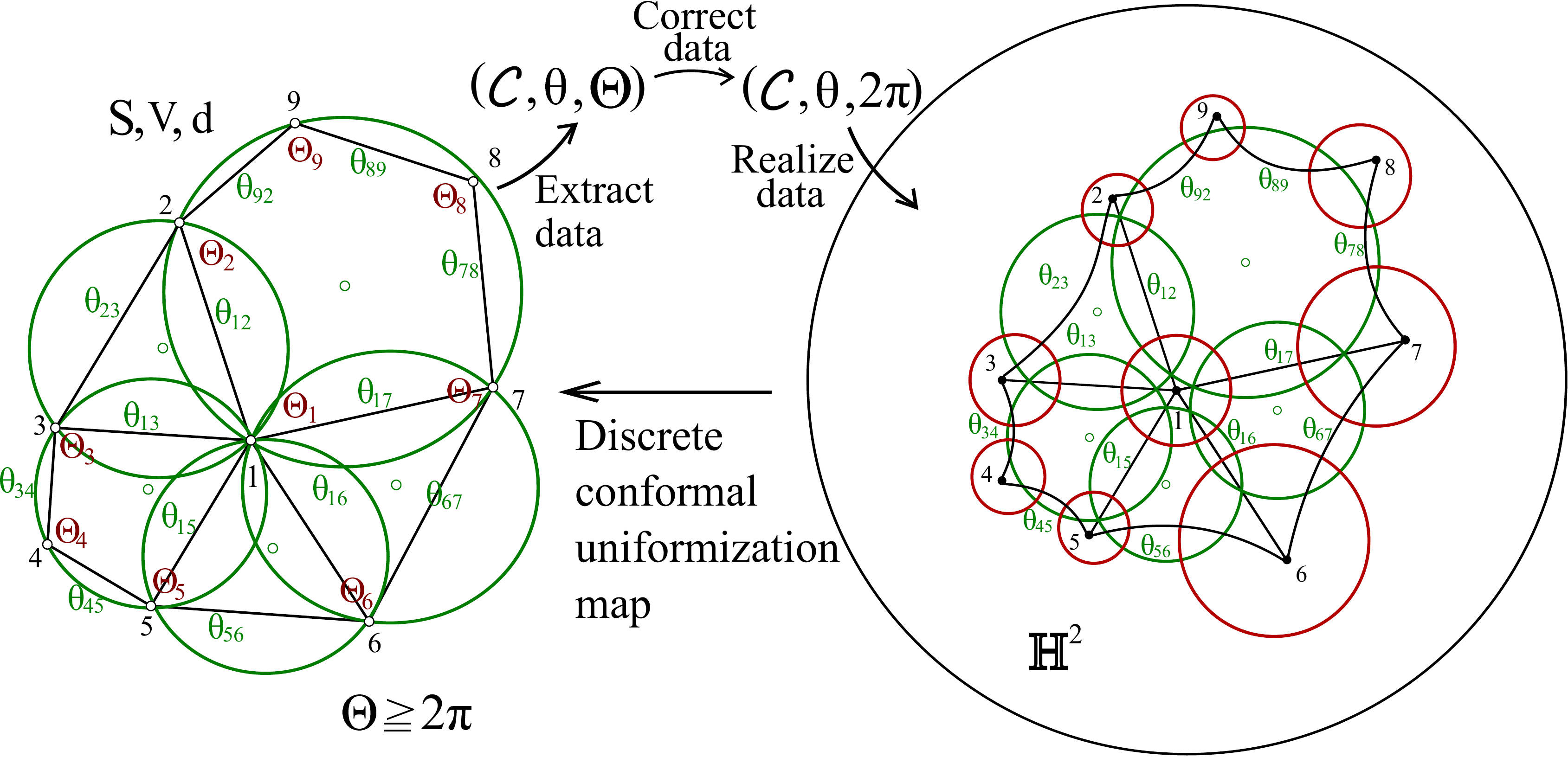}
\caption{Discrete uniformization according to Theorem
\ref{thm_uniformization_negative_curvature}} \label{Fig2}
\end{figure}

\begin{thm} \label{thm_uniformization_negative_curvature}
There exists a hyper-ideal circle pattern on $S$, with underlying
smooth hyperbolic metric $h$, that realizes the combinatorial
angle data $(\cellcomplex, \theta, 2\pi)$ (see Figure \ref{Fig2}).
The pattern is unique up to isometry isotopic to identity (i.e.
isometry preserving the combinatorics of the cell complex and its
labelling). In other words, the Delaunay circle pattern
corresponding to $(S,d,V)$ is discretely conformally equivalent to
a unique, up to label-preserving isometry, hyper-ideal circle
pattern with a smooth hyperbolic metric $h$ on $S$.
\end{thm}

Observe that Theorem \ref{thm_uniformization_negative_curvature}
implies the existence of a unique discrete conformal map (see
Definition \ref{Def_discrete_conformal_equivalence} and Figure
\ref{Fig2}) between the classical Delaunay circle pattern on
$(S,d,V)$ and a uniquely defined (up to isometry) hyper-circle
pattern with a smooth hyperbolic metric on $S$. Consequently, one
can develop the geodesic cell-complex corresponding to this newly
obtained hyper-ideal circle pattern in the hyperbolic plane and
thus obtain a fundamental domain whose isometric gluing defines a
Fuchsian group $\Gamma$ such that $\hyperbolicplane / \Gamma$ is
isometric to $S,h$. Moreover, one ends up with a discrete
conformal universal covering map from the $\Gamma$-invariant
hyper-ideal circle pattern on $\hyperbolicplane$ to the Delaunay
circle pattern on $S, d$ (Figure \ref{Fig2}). This is reminiscent
of the classical uniformization theorem for Riemann surfaces and
therefore we consider Theorem
\ref{thm_uniformization_negative_curvature} as a discrete version
of it.

\section{Discrete uniformization of branch covers over the Riemann sphere}
\label{Sec_uniformization_hyper_elliptic_surfaces}

In relation to surfaces with discrete non-positive curvature, we
propose a circle pattern uniformization of Riemann surfaces
represented as finite branch covers over the Riemann sphere.

Let $p \, : \, S \, \to \, \hat{\CC}$ be a finite topological
branch cover over the Riemann sphere $\hat{\CC}$. For example, one
can think of $S$ as a smooth compact algebraic curve in $\mathbb{CP}^2$ 
whose affine part is given by a complex polynomial equation
$P(x,y) = 0$.  The branch covering map $p$ can be defined as $p \,
: \, (x,y) \in S \, \mapsto \, x$ for example, although any
meromorphic function on $S$ would do.  
Define by $V_1$ all ramification points of $p$ on $S$ and let
$V_1(\hat{\CC}) = p(V_1)$ be the branch points of $p$ on
$\hat{\CC}$. Observe that $p^{-1}(V_1(\hat{\CC}))$ contains $V_1$
but does not necessarily coincide with it. Denote by $N \in
\naturals$ the number of sheets of $p$.

As input data we consider:

\smallskip

\noindent (a) a finite topological branch cover $p \, : \, S \,
\to \, \hat{\CC}$ with ramification points $V_1 \subset S$ and
branch points $V_1(\hat{\CC}) \subset \hat{\CC};$

\smallskip
\noindent (b) a finite set of points $V_{\hat{\CC}} =
V_0(\hat{\CC}) \cup V_1(\hat{\CC})$ on $\hat{\CC}$ where
$V_0(\hat{\CC}) \cap V_1(\hat{\CC}) = \varnothing$.

\smallskip
\noindent (c) a finite set $V = p^{-1}(V_{\hat{\CC}})$ on $S$ with
$V_0 = V \setminus V_1$.

\smallskip








\noindent This input data can be denoted by $ p : (S, V) \to
(\hat{\CC}, V_{\hat{\CC}})$. Observe that we do not need to
specify any metric on $\hat{\CC}$ because of its natural conformal
structure in which, as already mentioned in the Discussion on
$\hat{\CC}$ hyper-ideal patterns and the paragraph preceding it in
Section \ref{Sec_def_and_notations}, the notions of circles and
angles are invariant under the action of the conformal
automorphism group $\mathbb{P}SL(2,\CC)$, so Delaunay circle
patterns naturally exist on $\hat{\CC}$. Consequently, let us
generate the unique Delaunay circle pattern on $\hat{\CC}$ with
repsect to the points $V_{\hat{\CC}}$. Recall, that in the
Discussion of planar weighted Delaunay tessellations of Section
\ref{Sec_def_and_notations} we have explained how, after
stereographic projection, one can interpret hyper-ideal circle
patterns (and Delaunay ones in particular) from either Euclidean
or hyperbolic point of view, as specified in Definition
\ref{Def_hyperideal_circle_pattern}, and also obtain the
corresponding geodesic weighted Delaunay cell complexes and their
dual geodesic Voronoi cell complexes. Consequently, we have an
embedding of the topological complex $\cellcomplex_{\hat{\CC}} =
(V_{\hat{\CC}}, E_{\hat{\CC}}, F_{\hat{\CC}})$ and its dual
$\cellcomplex^*_{\hat{\CC}} = (V^*_{\hat{\CC}}, E^*_{\hat{\CC}},
F^*_{\hat{\CC}})$ in $\hat{\CC}$. By definition, a pair of
Delaunay circles (a special case of face circles) are adjacent to
each other exactly when there is a corresponding edge $ij \in
E_{\hat{\CC}}$ of $\cellcomplex_{\hat{\CC}}$, and equivalently a
dual edge $ij^* \in E_{\hat{\CC}}^*$ of
$\cellcomplex_{\hat{\CC}}^*$. Thus, the intersection angle between
them $\hat{\theta}_{ij} \in  (0,\pi)$ gives rise to an angle
assignment $\hat{\theta} \, : \, E_{\hat{\CC}} \, \to \, (0,\pi)$
(and by duality $\hat{\theta} \, : \, E_{\hat{\CC}}^* \, \to \,
(0,\pi)$). Since the branch points of $p : S \to \hat{\CC}$ are
among the vertices $V_{\hat{\CC}}$, we can lift
$\cellcomplex_{\hat{\CC}}$ and its dual
$\cellcomplex_{\hat{\CC}}^*$ to the embedded complexes
$\cellcomplex=(V,E,F)$ and $\cellcomplex^*=(V^*,E^*,F^*)$, which
are also dual to each other, i.e. $\cellcomplex =
p^{-1}(\cellcomplex_{\hat{\CC}})$ and $\cellcomplex^* =
p^{-1}(\cellcomplex_{\hat{\CC}}^*)$. The ramification points of
the covering are by construction among the vertices of
$\cellcomplex$ and thus lie in the interiors of their
corresponding dual faces from $\cellcomplex^*$. Furthermore, we
can lift the angle assignment $\hat{\theta}$ to the covering angle
assignment $\theta \, : \, E \, \to \, (0, \pi)$ by $\theta_{ij} =
\hat{\theta}_{p(ij)}$ for all $ij \in E$. By duality, we also have
$\theta \, : \, E^* \, \to \, (0, \pi)$.


With the preceding constructions in mind, one can extract from
$(\hat{\CC}, V_{\hat{\CC}})$ the combinatorial angle data
$(\cellcomplex, \theta, \Theta)$, where

\smallskip

\noindent 1. $\cellcomplex = (V,E,F)$ is the topological
cell-complex that branch covers the Delaunay cell complex
$\cellcomplex_{\hat{\CC}}$ via the covering map $p$, as described
above.

\smallskip
\noindent 2. $\theta \, : \, E \, \to \, (0,\pi)$  are the lifts
of the intersection angles between the pairs of adjacent Delaunay
circles, also defined above.


\smallskip
\noindent 3. $\Theta_k = 2\pi$ for all $k \in V_0$ and $\Theta_k =
2\pi N_k$ for all $k \in V_1$, where $N_k = 2,3,..., N$ is the
index of the ramification point $k$, i.e. the number of sheets
meeting at the ramification point.

\smallskip
To discretely uniformize, correct the data (See Figure \ref{Fig2})
to $(\cellcomplex, \theta, 2\pi)$.

\begin{thm} \label{thm_general_uniformization}
There exists a hyper-ideal circle pattern with an underlying
either (i) hyperbolic cone metric $h$ on $S$, when ${genus}(S)
\geq 2$, or (ii) a Euclidean cone metric $h$ on $S$, when
${genus}(S) = 1$, which realizes the combinatorial angle data
$(\cellcomplex, \theta, 2\pi)$. The pattern is unique up to
isometry isotopic to identity (i.e. isometry preserving the
combinatorics of the cell complex and its labelling) as well as
scaling in the Euclidean case. In other words, the lift of the
standard Delaunay circle pattern of $(\hat{\CC},V_{\hat{\CC}})$
via $p$ on $S$ is discretely conformally equivalent to a unique
(up to isometry isotopic to identity and scaling when Euclidean)
hyper-ideal circle pattern with either  a hyperbolic (genus $\geq
2$) or a Euclidean (genus $=1$) metric $h$ on $S$.
\end{thm}

We state separately the case of $S$ being homeomorphic to a
sphere.

\begin{thm} \label{thm_sphere_uniformization}
There exists a hyper-ideal circle pattern on $\hat{\CC}$, unique
up to label preserving $\mathbb{P}SL(2,\CC)$ transformation, which
realizes the combinatorial angle data $(\cellcomplex, \theta,
2\pi)$. In other words, the lift of the standard Delaunay circle
pattern of $(\hat{\CC},V_{\hat{\CC}})$ via $p$ on the topological
sphere $S$ is discretely conformally equivalent to a unique, up to
label respecting M\"obius transformation, hyper-ideal circle
pattern on the Reiamnn sphere $\hat{\CC}$.
\end{thm}

Theorems \ref{thm_general_uniformization} and
\ref{thm_sphere_uniformization} combined lead to a discrete
analogue of the uniformization theorem for closed Riemann surfaces
represented as branch covers over the Riemann sphere. In
particular, these two theorems can provide a discrete
uniformization method for smooth complex algebraic curves.



Just like in the case of the classical uniformization theorem,
when the genus of $S$ is two or greater, one can develop the final
hyper-ideal circle patterns from Theorems
\ref{thm_uniformization_negative_curvature} and
\ref{thm_general_uniformization} in the hyperbolic plane (e.g. the
upper-half plane model or the Poincar\'e disk model) and obtain a
Fuchsian group $\Gamma$, unique up to conjugation by a hyperbolic
isometry, as well as a $\Gamma-$invariant hyper-ideal circle
pattern on $\hyperbolicplane$ whose factor $\hyperbolicplane /
\Gamma$ is isometric to $(S,h)$ together with the hyper-ideal
circle pattern described in Theorems
\ref{thm_uniformization_negative_curvature} or
\ref{thm_general_uniformization}. When the genus of $S$ is one,
the same is true but for the Euclidean plane together with a
lattice group of translations on it.

\section{Existence and uniqueness of hyper-ideal patterns}
\label{Sec_Existence_Uniqueness_of_Patterns}

In this section we introduce the main tool used in the proof of
all theorems stated in the previous section.

As already discussed in the paragraph preceding Definition
\ref{Def_discrete_conformal_equivalence} from Section
\ref{Sec_def_and_notations}, given any hyper-ideal circle pattern
one can always extract from it the combinatorial data
$(\cellcomplex,\theta,\Theta)$, where $\cellcomplex$ is a
cell-complex representing the combinatorics of the pattern,
$\Theta : V \to (0, \infty)$ is the assignment of cone-angles at
the vertices of the complex and $\theta : E_d \to (0,\pi)$ is the
assignment of intersection angles between adjacent face circles of
the pattern. The proofs of the discrete uniformization theorems
heavily rely on the solution to the following problem: \emph{Find
a hyperbolic or flat cone-metric, together with a hyper-ideal
circle pattern on the given surface that realizes the data
$(\cellcomplex, \theta, \Theta)$.} 


The solution to this problem was first presented by Schlenker in
\cite{Sch1}. However, we will utilize the approach and the
notations used in \cite{ND}. We start with the constructions and
the definitions necessary for the formulation of the main result
in this section.

\begin{figure}
\centering
\includegraphics[width=14cm]{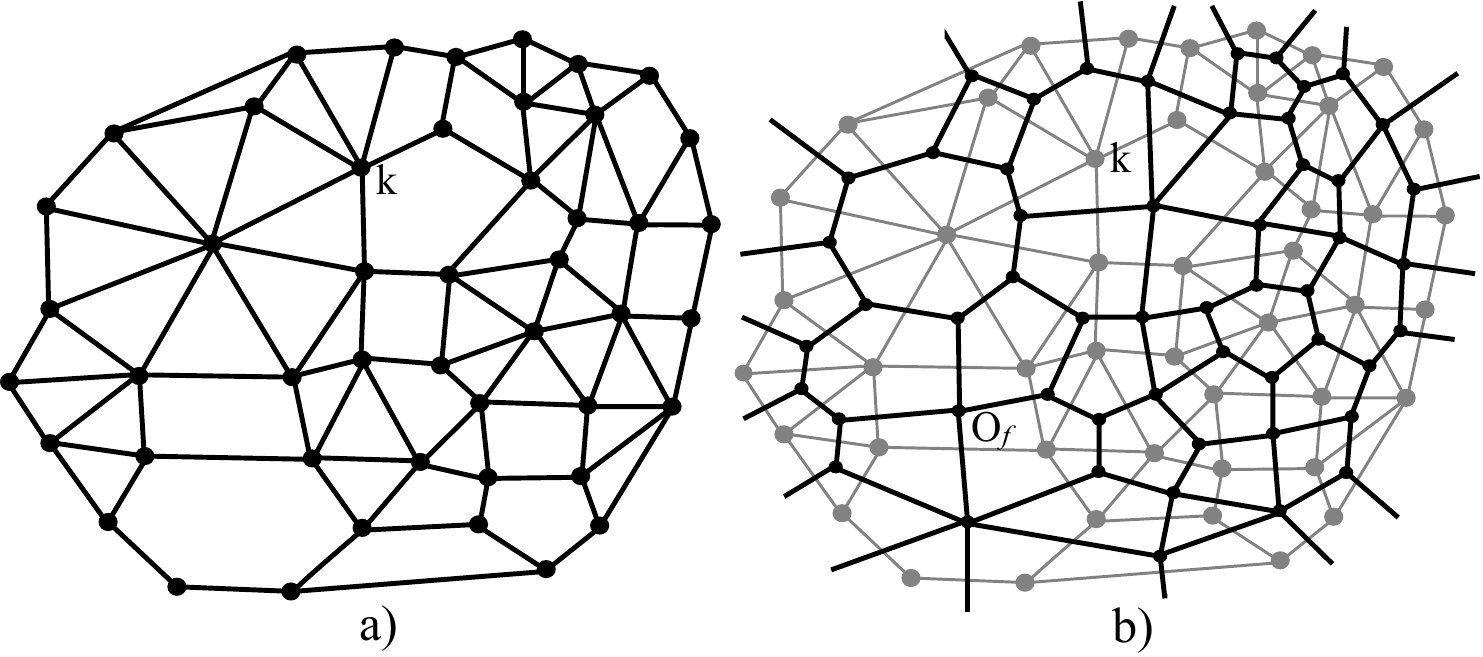}
\caption{a) The cell complex $\cellcomplex=(V,E,F)$ and b) its
dual $\cellcomplex^*=(V^*,E^*,F^*)$} \label{Fig3}
\end{figure}

Assume a cell complex $\cellcomplex = (V, E, F)$ is fixed on the
surface $S$ (Figure \ref{Fig3}a). Denote by $\cellcomplex^* =
(V^*, E^*, F^*)$ the cell complex dual to $\cellcomplex$, where
$V^*$ are the dual vertices, $E^*$ are the dual edges and $F^*$
are the dual faces (see Figure \ref{Fig3}b). On Figure \ref{Fig3}b
the elements of the original complex $\cellcomplex$ are drawn in
grey, while the elements of the dual complex $\cellcomplex^*$ are
in black.

\begin{figure}
\centering
\includegraphics[width=9.5cm]{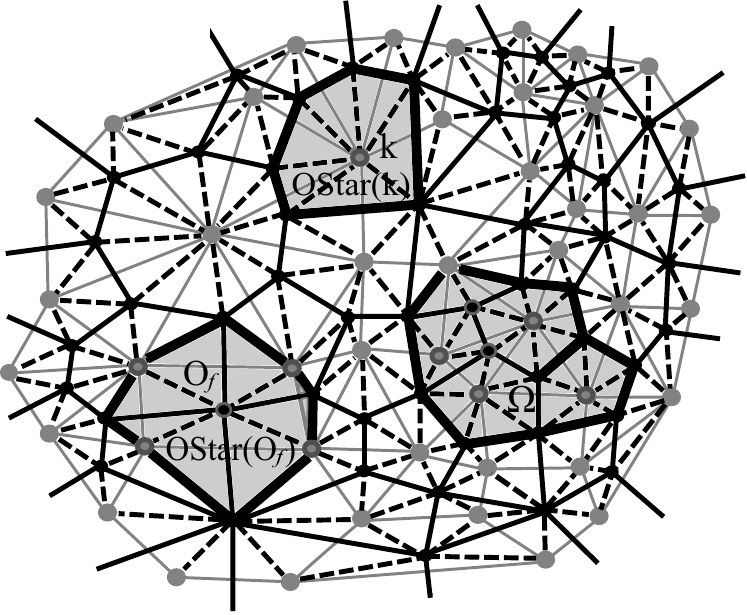}
\caption{The triangulation
$\hat{\Triang}=(\hat{V},\hat{E},\hat{F})$ together with two
examples of open stars and one admissible domain $\Omega$.}
\label{Fig4}
\end{figure}

Next, define the subdivision $\hat{\Triang}=(\hat{V}, \hat{E},
\hat{F})$ of $\cellcomplex^*$, shown on Figure \ref{Fig4}, where

\smallskip


\noindent $\bullet$ $\hat{V} = V \cup V^*$, i.e. the vertices of
$\hat{\Triang}$ consist of all vertices of $\cellcomplex$ and all
dual vertices. These are all black and grey vertices from Figure
\ref{Fig4};

\smallskip
\noindent $\bullet$ $\hat{E} = E^* \cup \big\{\, iO_f \,\, | \,\,
O_f \in V^* \,\, \text{and $i$ is a vertex of $f$ } \big\}$, i.e.
the edges of $\hat{\Triang}$ consist of all dual edges and all
edges, obtained by connecting a dual vertex $O_f \in f$ to all the
vertices of the face $f \in F$ it belongs to. The latter type of
edges will be called \emph{corner edges}. The dual edges can be
seen on both Figures \ref{Fig3}b and \ref{Fig4} painted solid
black, while the corner edges are the black dashed edges from
Figure \ref{Fig4}.

\smallskip
\noindent $\bullet$ $\hat{F} = \bigl\{ iO_fO_{f'} \, | \, \text{
$ij \in E$ common edge for $f$ and $f'$ from $F$ } \bigr\},$ i.e.
the faces of $\hat{\Triang}$ are the topological triangles
obtained by looking at the connected components of the complement
of the topological graph $(\hat{V}, \hat{E})$ on $S$. On Figure
\ref{Fig4} these are the triangles with one solid black and two
dashed black edges. They also have two black (dual) vertices and
one grey vertex.

The next important notion to be defined is the \emph{open star} of
a vertex from $\hat{\Triang}$.

\begin{Def} \label{Def_open_star}
Let $\hat{v} \in \hat{V}$ be an arbitrary vertex of
$\hat{\Triang}.$ Then its \emph{open star}
$\text{\emph{OStar}}(\hat{v})$ is defined as the open interior of
the union of all closed triangles from $\hat{\Triang}$ which
contain $\hat{v}$.
\end{Def} 

Whenever $\hat{v} = k \in V$ is a vertex of $\cellcomplex$, then
its open star is simply the open interior of the face from
$\cellcomplex^*$
dual to $k$. An example denoted by $\text{OStar}(k)$ and colored
in grey is shown on Figure \ref{Fig4}. The boundary of
$\text{OStar}(k)$ consists entirely of dual edges from $E^*$. If
we denote by $E_k$ the set of all edges of $\cellcomplex$ which
have vertex $k$ as an endpoint, then $\partial \, \text{OStar}(k)
= \cup \bigl\{ ik^* \in E^* \,\, | \,\, ik \in E_k \, \bigr\}$. If
$\hat{v} = O_f \in V^*$ is a vertex from the dual complex
$\cellcomplex^*$, then the boundary of its open star consists
entirely of corner edges from $\hat{\Triang}$ (e.g. the grey
region $\text{OStar}(O_f)$ on Figure \ref{Fig4}).

Recall that the vertex set $V$ of the cell-complex $\cellcomplex$
is always partitioned into two subsets $V_0$ and $V_1$.  


Following the terminology of \cite{Sch1} (see also \cite{ND}), one
can define what Schlenker calls an admissible domain.

\begin{Def} \label{Def_admissible_domain}
An open connected subdomain $\Omega$ of the surfaces $S$ is called
an \emph{admissible domain of $(S,\cellcomplex)$ } whenever the
following conditions hold:

\smallskip

\noindent {\bf 1.} There exists a subset $\hat{V}_0 \subseteq
\hat{V}$, such that $\Omega = \cup \big\{
\mbox{\emph{OStar}}(\hat{v}) \,\, | \,\, \hat{v} \in \hat{V}_0
\big\};$

\smallskip

\noindent {\bf 2.} $\Omega \neq \varnothing$ and $\Omega \neq S$
and $\Omega \cap V \neq \varnothing$.

\smallskip

\noindent {\bf 3.} $\Omega$ is not punctured, i.e. if
$\emph{OStar}(\hat{v}) \setminus \{\hat{v}\} \subseteq \Omega$
then $\hat{v} \in \Omega$ for any $\hat{v} \in \hat{V}$.
Consequently, the boundary of $\Omega$ is a nonempty set of edges
of $\hat{\Triang}$, i.e. dual and/or corner edges.
\end{Def}

A special example of an admissible domain is the open star of a
vertex of $\cellcomplex$. The open star of a dual vertex however
is not an admissible domain because it is disjoint from $V$. An
example of an admissible domain can be seen on Figure \ref{Fig4},
denoted by the symbol $\Omega$ and shaded in grey. On this picture
$\Omega$ is simply connected but in general it doesn't have to be.

The boundary of an admissible domain $\Omega$ is a disjoint union
of topological curves on $S$, consisting entirely of edges
together with their vertices, belonging to the triangulation
$\hat{\Triang}$. In other words, the boundary of $\Omega$ consists
of dual edges and/or corner edges from $\hat{E}$ (and their
vertices), but all of its connected components are interpreted as
closed paths in the one-skeleton of $\hat\Triang$, so that some of
the edges could be traced (counted) twice (see Figure \ref{Fig4})
as well as some of the vertices could be counted twice or more
times. Edges are traced twice exactly when an edge of
$\hat{\Triang}$ is disjoint from $\Omega$, but the interiors of
the two topological triangles from $\hat{\Triang}$, lying on both
sides of the edge, are contained in $\Omega$. We denote this
version of the boundary of $\Omega$ by $\partial\Omega$.

\begin{thm} \label{Thm_main}\emph{(Schlenker \cite{Sch1,ND})}
Let $S$ be a closed surface with a topological cell complex
$\cellcomplex = (V, E, F)$ on it. Assume also that $V = V_0 \sqcup
V_1$. Then the combinatorial angle data $(\cellcomplex, \theta,
\Theta)$ is realized by a hyperbolic or Euclidean hyper-ideal
circle pattern on $S$ if and only if the assignment of angles
$(\theta, \Theta) \in \reals^E \times \reals^V$ satisfies the
following conditions:


\medskip

\noindent 1) $\theta_{ij} \in (0,\pi)$ for any $ij \in E$;

\medskip

\noindent 2) $\Theta_k > 0$ for all $k \in V_1$ and $\Theta_k =
\sum_{ik \in E_k} (\pi - \theta_{ik})$ for all $k \in V_0$. The
latter can be also written as $\Theta_k = \sum_{ik^* \subset
\partial\Omega} (\pi - \theta_{ik})$ for $\Omega =
\mbox{\emph{OStar}}(k)$;

\medskip

\noindent 3) $\sum_{k \in V} (2 \pi - \Theta_k) > 2 \pi \chi(S)$
in the hyperbolic case and $\sum_{k \in V} (2 \pi - \Theta_k) = 2
\pi \chi(S)$ in the Euclidean case;

\medskip

\noindent 4) For any admissible domain $\Omega$ of $(S,
\cellcomplex)$, such that $\Omega \neq \mbox{\emph{OStar}}(k)$ for
some $k \in V_0$,
\begin{equation} \label{Eqn_first_condition}
\sum_{ij^* \subset \partial\Omega} (\pi - \theta_{ij}) + \sum_{k
\in \Omega \cap V} (2 \pi - \Theta_k) + \pi |\partial\Omega \cap
V| > 2 \pi \chi(\Omega).
\end{equation}
The notations $\chi(S)$ and $\chi(\Omega)$ are the Euler
characteristics of $S$ and $\Omega$ respectively.

\medskip
\noindent Furthermore, whenever the hyper-ideal circle pattern in
question exists, it is unique up to isometry isotopic to identity,
as well as scaling in the Euclidean case. Finally, the hyper-ideal
circle pattern can be reconstructed from the unique critical point
of a strictly convex functional defined on an open convex set,
bounded only by coordinate hyperplanes of $\reals^N$ for some
suitable $N \in \naturals$.
\end{thm}

Conditions 1 to 4 from Theorem \ref{Thm_main} describe a convex
polytope, which in the hyperbolic case we denote by
$\polytope^{h}$ and in the Euclidean case by $\polytope^{e}$. We
will also use the common notation $\polytope$ whenever we do not
need to specify the underlying geometry. We call these polytopes
\emph{angle data polytopes}.

To optimize the conditions of Theorem \ref{Thm_main} a bit more,
one can define the so called \emph{strict admissible domain}.

\begin{Def} \label{Def_strict_admissible_domain}
An open connected subdomain $\Omega$ on the surfaces $S$ is called
a \emph{strict admissible domain of $(S,\cellcomplex)$ } whenever
$\Omega$ is admissible and $\partial\Omega \cap V_0 =
\varnothing$.
\end{Def}

As shown in \cite{ND}, the angle data polytopes can be described
via strict admissible domains instead of admissible domains.
Simply, the admissible domains which are not strict do not add
more restrictions to the angle data.

\begin{cor} \label{Cor_polytopes_general_admissible_domains}
The statements of Theorem \ref{Thm_main} still hold even when the
expression ``admissible domain" in point 4) of Theorem
\ref{Thm_main} is replaced by the expression ``strict admissible
domain".
\end{cor}

\section{Proof of Theorem \ref{thm_uniformization_negative_curvature}}
\label{Sec_proof_uniformization_negative_curvature} 

Let $(\cellcomplex, \theta, \Theta)$ be the combinatorial angle
data extracted from the Delaunay circle pattern on $S,d$ with
respect to the finite set of points $V$ (refer to the explanations
towards the end of Section \ref{Sec_def_and_notations}. See also
Figure \ref{Fig2}). This circle pattern is a special case of a
hyper-ideal circle pattern with all vertex circles of radius zero.
Hence, Theorem \ref{Thm_main} applies to it. Then the angles
$(\theta, \Theta)$ satisfy the conditions

\noindent(i) $\theta_{ij} \in (0,\pi)$ for any $ij \in E$;

\smallskip

\noindent (ii) $\Theta_k = \sum_{ik \in E_k} (\pi - \theta_{ik})
\geq 2 \pi$ for any $k \in V$. Recall that $V$ has been split into
$V_0$ containing all $k$ such that $\Theta_k = 2\pi$, and $V_1$
consisting of all $k$ for which $\Theta_k > 2\pi$. By assumption,
$V = V_0 \cup V_1$.

\smallskip

\noindent (iii) $\sum_{k \in V} (2 \pi - \Theta_k) \geq 2 \pi
\chi(S)$, where in the case of Euclidean cone-metric $d$ we have
an equality, and in the case of a hyperbolic cone-metric we have a
strict inequality.

\smallskip

\noindent (iv) For any admissible domain $\Omega$ of $(S,
\cellcomplex)$, such that $\Omega \neq \text{OStar}(k)$ for some
$k \in V_0$,

\begin{equation} \label{Eqn_first_condition_revisited}
\sum_{ij^* \subset \partial\Omega} (\pi - \theta_{ij}) + \sum_{k
\in \Omega \cap V} (2 \pi - \Theta_k) + \pi |\partial\Omega \cap
V|
> 2 \pi \chi(\Omega).
\end{equation}

\smallskip

To prove Theorem \ref{thm_uniformization_negative_curvature} all
we have to do is check whether the angle data $(\theta, 2\pi)$
belongs to the polytope $\polytope^h$, whose description is given
by points 1, 2, 3 and 4 from Theorem \ref{Thm_main}. If we manage
to confirm that, Theorem \ref{Thm_main} guarantees the existence
and the uniqueness of the sought hyper-ideal circle pattern with a
smooth underlying hyperbolic metric $h$ on $S$. So let us check
conditions 1 to 4 of Theorem \ref{Thm_main}. We point out that if
the metric $d$ is already smooth, then we have nothing to prove,
as $d$ is actually the metric $h$ and the hyper-ideal circle
pattern in question is in fact the original Delaunay circle
pattern. Therefore, for the rest of this proof we assume that $d$
has at least one cone-singularity with cone angle greater than
$2\pi$, which means that $V_1 \neq \varnothing$.

Since the data $(\cellcomplex, \theta, \Theta)$ comes from a
Delaunay circle pattern, it is clear that condition 1 is
immediately true. Condition 2 is also true because if $k \in V_0$,
then $\Theta_k = 2\pi = \sum_{ik \in E_k} (\pi - \theta_{ik})$
follows from (ii) above. The case when $k \in V_1$ implies that
$\Theta_k
> 2 \pi$ and therefore the already established in (ii) above
strict inequality $\sum_{ik \in E_k} (\pi - \theta_{ik}) > 2\pi$
can be rewritten as $\sum_{ik^* \in \partial \Omega} (\pi -
\theta_{ik}) + (2\pi - 2\pi)
> 2\pi = 2\pi \chi(\Omega)$ which is a special case of condition 4 with
$\Omega = \text{OStar}(k)$.

Since $V_1 \neq \varnothing$ and since we have assumed in Theorem
\ref{thm_uniformization_negative_curvature} that $\Theta_k \geq
2\pi$ holds for each $k \in V$, then there is at least one $k \in
V$ for which $\Theta_k
> 2\pi$. Therefore $0 = \sum_{k \in V} (2\pi - 2\pi) > \sum_{k \in
V} (2\pi - \Theta_k) \geq 2\pi \chi(S)$ by (iii) above. Thus,
condition 3 is confirmed.

Finally, in order to verify point 4) from Theorem \ref{Thm_main},
let us fix an arbitrary admissible domain $\Omega$ (see Definition
\ref{Def_admissible_domain}). Using again the assumption that
$\Theta_k \geq 2\pi$, for all $k \in V$, 
and combining it with observation (iv) above, we can deduce that

\begin{equation*}
\sum_{ij^* \subset \partial\Omega} (\pi - \theta_{ij}) + \sum_{k
\in \Omega \cap V} (2 \pi - 2 \pi) \geq \sum_{ij^* \subset
\partial\Omega} (\pi - \theta_{ij}) + \sum_{k \in \Omega \cap V}
(2 \pi - \Theta_k) > 2 \pi \chi(\Omega) - \pi |\partial\Omega \cap
V|.
\end{equation*}

The proof of Theorem \ref{thm_uniformization_negative_curvature}
is complete.


\section{Realization theorems for circle patterns on the sphere}
\label{Sec_realization_patterns_on_sphere}

One of the central tools in the proof of Theorems
\ref{thm_general_uniformization} and
\ref{thm_sphere_uniformization}, alongside Theorem \ref{Thm_main},
are the Theorem of Rivin \cite{Riv96} and the Theorem of Bao and
Bonahon \cite{BaoBon}. Originally, both results are stated as a
geometric characterization of (i) convex ideal polyhedra (in the
case of Rivin) and (ii) convex hyper-ideal polyhedra (in the case
of Bao and Bonahon) in the hyperbolic three-space. However, due to
the natural one-to-one correspondence between such polyhedra and
(i) Delaunay circle patterns (Rivin \cite{Riv96}) and (ii)
hype-ideal circle patterns (Bao and Bonahon \cite{BaoBon}) on the
Reimann sphere (see also \cite{ThuNotes,BobSpr1,Schl} as well as
Section \ref{Sec_Spaces_of_triangles_and_patterns}), we restate
them as theorems about circle patterns.

We begin with some notations and definition. Let $\cellcomplex =
(V,E,F)$ be a topological cell complex and let $\cellcomplex^* =
(V^*, E^*, F^*)$ be its dual on a surface $S$ (Figure \ref{Fig3}).
Denote by $\cellcomplex^*(1)$ the one-skeleton of
$\cellcomplex^*$, which is the subcomplex of $\cellcomplex^*$
composed only of its vertices and edges (so faces excluded). In
other words, $\cellcomplex^*(1)=(V^*,E^*)$ and so it is a graph
embedded in $S$. A \emph{path in the one-skeleton} of
$\cellcomplex^*$ is a sequence of dual edges from $E^*$ whose
union is a topological path on $S$. Analogously, a \emph{loop in
the one-skeleton} of $\cellcomplex^*$ is a sequence of dual edges
from $E^*$ whose union is a topological loop on $S$. A path or a
loop in $\cellcomplex^*(1)$ is called \emph{simple} if it is
respectively a simple topological path or a simple loop on $S$.
This means that no two edges in the sequence repeat and each dual
vertex that the simple path or loop goes through is adjacent to no
more than two edges from the path (exactly two in the case of a
simple loop). The same terminology applies to the triangulation
$\hat{\Triang}$ (Figure \ref{Fig4}).

In the case when $S$ is a topological sphere we have the following
two powerful theorems.

\begin{thm} \label{thm_Rivin}
\emph{(Rivin \cite{Riv96})} Let $S$ be a topological sphere with a
strongly regular cell complex $\cellcomplex=(V,E,F)$ on it. Let
$\cellcomplex^*$ be the dual complex of $\cellcomplex$. Then the
combinatorial data $(\cellcomplex, \theta, 2\pi)$ is realized by a
Delaunay circle pattern on $\hat{\CC}$ if and only if the
following conditions are satisfied by the angle assignment $\theta
\, : \, E \, \to \, \reals$:

\smallskip


\noindent 1) $\theta_{ij} \in (0,\pi)$ for any $ij \in E$;

\smallskip

\noindent 2) For each simple loop $\delta$ in the one-skeleton of
$\cellcomplex^*$
$$ \sum_{ij^* \subset \delta} (\pi - \theta_{ij}) \geq 2 \pi$$
where equality holds if and only if $\delta$ is the boundary of a
face from $\cellcomplex^*$.


\smallskip

\noindent Furthermore, whenever the hyper-ideal circle pattern in
question exists, it is unique up to a $\mathbb{P}SL(2,\CC)$
transformation respecting the labelling of the complex.
\end{thm}

\noindent Rinin's Theorem is a special case of the more general
Bao and Bonahon's Theorem.

\begin{thm} \label{thm_Bao}
\emph{(Bao and Bonahon \cite{BaoBon})} Let $S$ be a topological
sphere with a strongly regular cell complex $\cellcomplex=(V,E,F)$
on it. Assume the vertices are partitioned into $V=V_0 \sqcup
V_1$. Furthermore, let $\cellcomplex^*$ be the dual complex of
$\cellcomplex$. Then the combinatorial data $(\cellcomplex,
\theta, 2\pi)$ is realized by a hyper-ideal circle pattern on
$\hat{\CC}$ with true vertex circles corresponding to the vertices
from $V_1$ and circles shrunk to points corresponding to the
vertices from $V_0$ if and only if the following conditions are
satisfied by the angle assignment $\theta \, : \, E \, \to \,
\reals$:

\smallskip


\noindent 1) $\theta_{ij} \in (0,\pi)$ for any $ij \in E$;

\smallskip

\noindent 2) For each simple loop $\delta$ in the one-skeleton of
$\cellcomplex^*$
$$ \sum_{ij^* \subset \delta} (\pi - \theta_{ij}) \geq 2 \pi$$
where equality holds if and only if $\delta$ is the boundary of a
face from $\cellcomplex^*$ dual to a vertex from $V_0$.


\smallskip

\noindent 3) For each simple path $\delta$ in the one-skeleton of
$\cellcomplex^*$ joining two distinct dual vertices of the same
dual face so that $\delta$ is not contained in the boundary of any
dual face from $\cellcomplex^*$
$$ \sum_{ij^* \subset \delta} (\pi - \theta_{ij}) > \pi.$$

\smallskip

\noindent Furthermore, whenever the hyper-ideal circle pattern in
question exists, it is unique up to a $\mathbb{P}SL(2,\CC)$
transformation respecting the labelling of the complex.
\end{thm}


\section{Proof of Theorem \ref{thm_general_uniformization}}
\label{Sec_proof_uniformization_branch_covers}


Let $\hat{\Triang}_{\hat{\CC}} = (\hat{V}_{\hat{\CC}},
\hat{E}_{\hat{\CC}}, \hat{F}_{\hat{\CC}})$ be the triangular
subdivision of the dual complex $\cellcomplex^*$ as defined in
Section \ref{Sec_Existence_Uniqueness_of_Patterns} and shown on
Figure \ref{Fig4}. Lift $\hat{\Triang}_{\hat{\CC}}$ to the
triangulation $\hat{\Triang}$ on the surface $S$ via the branch
covering map $p$. By construction, $\hat{\Triang}$ is the
subtriangulation of the lifted dual complex $\cellcomplex^*$ on
$S$ described in Section
\ref{Sec_Existence_Uniqueness_of_Patterns} and depicted on Figure
\ref{Fig4}. Thus, $p(\hat{\Triang}) = \hat{\Triang}_{\hat{\CC}}$
and $p(\cellcomplex^*) = \cellcomplex^*_{\hat{\CC}}$. Just like
before, in order to prove Theorem \ref{thm_general_uniformization}
we simply have to check whether the angle data $(\theta, 2\pi)$
belongs to the polytope $\polytope$, whose description is given by
points 1, 2, 3 and 4 of Theorem \ref{Thm_main}. However, in this
case we are going to apply Theorem \ref{Thm_main} in the setting
of Corollary \ref{Cor_polytopes_general_admissible_domains}, which
means we are going to use the description of the polytope
$\polytope$ in terms of strict admissible domains instead of
admissible ones.

\subsection{Verification of conditions 1, 2 and 3 of
Theorem \ref{Thm_main}}

\paragraph{Verification of conditions 1 and 3.} Since the angles $\theta = p^*\hat{\theta} \in \reals^{E}$ are
lifts of the angles of a Delaunay circle pattern on $\hat{\CC}$,
condition 1 of Theorem \ref{Thm_main} is automatically satisfied.
Denoting by $g(S)$ the genus of the surface $S$, if we assume that
$g(S) \geq 1$ then condition 3 is also satisfied because we have
taken $\Theta_k = 2\pi$ for all $k \in V$ and thus
$$0 \, = \, \sum_{k \in V} (2\pi - 2\pi) \, \geq \, 2\pi \chi(S) \, = \,
4\pi (1-g(S)).$$ The inequality is strict when $S$ has genus at
least two and becomes an equality when $S$ has genus 1, i.e. the
case of the torus.

\paragraph{Verification of condition 2.} From now on, whenever $k$ is a vertex of $\cellcomplex$ or
$\cellcomplex_{\hat{\CC}}$, we denote by $\text{Star(k)}$  the
closure of its open star $\text{OStar}(k)$ (for the definition of
$\text{OStar}(k)$ see Section
\ref{Sec_Existence_Uniqueness_of_Patterns} and Figure \ref{Fig4}).
This simply means that we add to the open star all dual edges
lying on its boundary. By definition, it is equivalent to saying
that the open star of $k$ is the interior of the face dual to $k$
and so its closure is the closure of that dual face. Since we have
assumed that the complexes are strongly regular, $\text{Star}(k)$
is an embedded closed disk.

Recall that since $p$ is a branch covering map, the map $p \, : \,
S \setminus p^{-1}\big(V_1({\hat{\CC}})\big) \, \to \, \hat{\CC}
\setminus V_1({\hat{\CC}})$ defines a regular cover. Assume that
$k \in V_0$. First, if $k \in V_0 \setminus
p^{-1}\big(V_1({\hat{\CC}})\big)$ then $k$ lies on the regular
cover and by construction the restricted map $p|_{\text{Star}(k)}
\, : \,  \text{Star}(k) \, \to \, \text{Star}(p(k))$ is a
homeomorphism. Second, if $k \in V_0 \cap
p^{-1}\big(V_1({\hat{\CC}})\big)$ then $p|_{\text{Star}(k)} \, :
\,  \text{Star}(k) \, \to \, \text{Star}(p(k))$ is a branch
covering map between closed disks with only one ramification point
$k$ in the domain's interior with ramification index $N_k = 1,2,3,
..., N$. However, since $k$ is not an actual branch point, the
index is $N_k=1$ and so the restricted map $p|_{\text{Star}(k)}$
is again a homeomorphism. Consequently, since the boundaries
$\partial \, \text{OStar}(k)$ and $\partial\, \text{OStar}(p(k))$
are homeomorphic with the same number of dual edges, and also
$\theta_{ij} = \hat{\theta}_{p(ij)}$,
$$\sum_{ij^* \subset \partial
\text{OStar}(k)} (\pi - \theta_{ij}) \, = \, \sum_{uv^* \subset
\partial \text{OStar}(p(k))} (\pi - \hat{\theta}_{uv})\,  = \,
2\pi,$$ where the last equality follows from condition 2 of
Rivin's Theorem \ref{thm_Rivin} (as well as condition 2 of Bao and
Bonahon's Theorem \ref{thm_Bao}). Thus, condition 2 of Schlenker's
Theorem \ref{Thm_main} is verified.



\subsection{Verification of condition 4 of Theorem \ref{Thm_main}}

\subsubsection{The case of admissible domains with general topology.}
Next, we focus on the verification of condition 4 of Theorem
\ref{Thm_main}. Assume $\Omega$ is a strict admissible domain of
$S, \cellcomplex$ which is not the open star of a vertex from
$V_0$, bit it can be the open star of a vertex from $V_1$. By
Definitions \ref{Def_admissible_domain} and
\ref{Def_strict_admissible_domain}, $\Omega$ is homeomorphic to
the interior of a compact surface with boundary. Therefore its
Euler characteristic is $\chi(\Omega) = 2 - 2 H - B$, where $H$ is
the number of handles and $B$ is the number of boundary components
of $\Omega$. Observe that $B \geq 1$ because $\partial \Omega \neq
\varnothing$. Consequently, $\chi(\Omega) = 2 - 2H - B$ yields the
restriction $\chi(\Omega) \leq 1$ where $\chi(\Omega) = 1$ exactly
when $\Omega$ has one boundary component and no handles, which
means that it is an open topological disk. In all other cases
$\chi(\Omega) \leq 0$. Since $\partial \Omega \neq \varnothing$,
either $\partial \Omega \cap E^* \neq \varnothing$ or $\partial
\Omega \cap V_1 \neq \varnothing$ or both. Therefore whenever
$\chi(\Omega) < 1$,
$$\sum_{ij^* \subset \partial \Omega} (\pi - \theta_{ij}) \, + \,
\pi |\partial \Omega \cap V_1| \, >  \, 0 \, \geq \, 2 \pi
\chi(\Omega).$$ The latter inequality is condition 4 of Theorem
\ref{Thm_main} for strict admissible domains that are not open
topological disks. Hence, we focus on the case when $\Omega$ is an
open topological disk and so $\chi(\Omega) = 1$.

\subsubsection{The case of admissible topological disks.}
Denote by $\gamma = \partial \Omega$ the loop in the one-skeleton
of $\hat{\Triang}$ that traverses the boundary of the topological
disk $\Omega$. Recall that we interpret $\gamma=(\hat{e}_1,
\hat{e}_2,..., \hat{e}_n, \hat{e}_1)$ as a cyclic sequence of
edges of $\hat{\Triang}$, so an edge $\hat{e}_s \in \hat{E}$ could
appear twice or the loop may pass through the same vertex several
times. Let us denote by $|\gamma|$ the one dimensional subcomplex
of the one-skeleton of $\hat{\Triang}$ formed by the edges and
vertices on $\gamma$, so repetitions of edges and vertices are
ignored. Hence, while $\gamma$ is a sequence of edges (i.e.
ordered), $|\gamma|$ is a one-complex with vertices and edges
(unordered). Absolutely the same notation we use in the case when
$\gamma$ is a path in $\hat{\Triang}(1)$.

The image of $\gamma$ under the branch covering map $p$ will be
denoted by $\delta = p(\gamma)$ and it will be interpreted as
$\delta=\big(p(\hat{e}_1), p(\hat{e}_2),..., p(\hat{e}_n),
p(\hat{e}_1)\big)$ which is a loop in the one-skeleton of
$\hat{\Triang}_{\hat{\CC}}$. Notice that the number of edges in
$\gamma$ and $\delta$ (counting repetitions) is the same.

\paragraph{Case 1.} Let $\partial \Omega \cap V_1 = \varnothing$.
Then $\gamma = \partial \Omega$ consists entirely of dual edges,
i.e. it is a path in the one-skeleton of the dual complex
$\cellcomplex^*$. Consequently, $\gamma$ lies on the regular cover
$p \, : \, S \setminus p^{-1}\big(V_1({\hat{\CC}})\big) \, \to \,
\hat{\CC} \setminus V_1({\hat{\CC}})$. Furthermore, its image
$\delta=p(\gamma)$ is a loop in the one-skeleton of the dual
complex $\cellcomplex_{\hat{\CC}}^*$ on the Riemann sphere.

\medskip

{\bf Subcase 1.1.} Let $|\delta|$ be a simple closed curve. Recall
the difference between $\delta$ and $|\delta|$. While the
one-dimensional subcomplex $|\delta|$ of the one-skeleton of
$\cellcomplex^*_{\hat{\CC}}$ is a simple closed curve on
$\hat{\CC}$, the sequence of edges $\delta$ is not necessarily
simple and it may traverse $|\delta|$ several times. However, the
assumption that $|\delta|$ is simple implies that there exists a
simple closed loop $\delta'$ in the one skeleton of
$\cellcomplex^*$ (i.e. a cyclic subsequence of $\delta$ without
edge repetitions), such that $|\delta'|=|\delta|$. Then $\hat{\CC}
\setminus |\delta| = \Omega_1 \cup \Omega_2$ where $\Omega_1$ and
$\Omega_2$ are disjoint open topological discs on the Riemann
sphere such that $\partial \Omega_1 = \partial \Omega_2 =
\delta'$.

\medskip

{\bf Subcase 1.1.A.} Assume that neither $\Omega_1$ nor $\Omega_2$
is an open star of a vertex from $V_{\hat{\CC}}$. Then by
condition 2 of Theorem \ref{thm_Rivin}, or alternatively Theorem
\ref{thm_Bao},
$$\sum_{uv^* \in \delta'} (\pi - \hat{\theta}_{uv})  > 2\pi.$$
 Since the map $p$ is onto and
$\delta=p(\gamma)$ as well as $\theta_{ij} = \hat{\theta}_{p(ij)}$
for all $ij \in E$, it is immediate to conclude that
$$\sum_{ij^* \in \gamma} (\pi - \theta_{ij}) \, = \,
\sum_{uv^* \in \delta} (\pi - \hat{\theta}_{uv}) \, \geq \,
\sum_{uv^* \in \delta'} (\pi - \hat{\theta}_{uv}) \, > \, 2\pi.$$

\medskip

{\bf Subcase 1.1.B.} We claim that under the assumptions of
condition 4 of Theorem \ref{Thm_main}, neither $\Omega_1$ nor
$\Omega_2$ can be the open star of a vertex from $V_0(\hat{\CC})$.
Indeed, assume that one of the two domains, say $\Omega_2$, is an
open star $\text{OStar}(\tilde{k}_0)$, where $\tilde{k}_0 \in
V_0(\hat{\CC})$. Then its closure is $\overline{\Omega}_2 =
\text{Star}(\tilde{k}_0)$ and it lies in the target space of a
regular cover. Consequently, due to the contractibility of
$\text{Star}(\tilde{k}_0)$ and the lifting property of covering
spaces, the preimage $p^{-1}(\text{Star}(\tilde{k}_0))$ is a
disjoint union of closed stars $\text{Star}({k}_s)$ for $s =
1...N$, where $k_s \in V_0$. The restriction of $p$ on each
$\text{Star}({k}_s)$ is a homeomorphism. Consequently, the full
preimage of $\delta$ is the disjoint union of the boundaries
$\partial \text{Star}(k_s)$ for $s=1...N$. Since
$\delta=p(\gamma)$, where $\gamma$ is a connected closed loop in
$\cellcomplex^*(1)$ which bounds an admissible domain on $S$,
$\gamma$ must be among the loops $\partial \text{Star}(k_s)$, i.e.
$\gamma = \partial \text{Star}(k_t)$ for some specific $t \in \{1,
..., N\}$ with $k_t \in V_0$. Consequently the loop $\gamma$ is
simple and it splits the surface $S$ into two open subdomains
$\text{OStar}(k_t)$ and $S \setminus \text{Star}(k_t)$. On the one
hand, we have assumed that the admissible domain $\Omega$ is not
the open star of a vertex from $V_0$ which means that $\Omega$ is
not $\text{OStar}(k_t)$. On the other hand, we have assumed that
$\Omega$ is a topological disc, while $S \setminus
\text{Star}(k_t)$ is definitely not a disc but a surface of genus
at least one with one closed disk removed (so it has at least one
handle). Thus, we have arrived at a contradiction due to the
assumption that $\Omega_2$ is the open star of $\tilde{k}_0 \in
V_0(\hat{\CC})$. Hence, this situation cannot occur.

\medskip

{\bf Subcase 1.1.C.} One of the two open discs, say $\Omega_1$, is
the open star $\text{OStar}(\tilde{k}_1)$ of a vertex $\tilde{k}_1
\in V_1(\hat{\CC})$. Then, by construction of the complexes
$\cellcomplex^*$ and $\hat{\Triang}$, the preimage of the closed
star $p^{-1}(\text{Star}(\tilde{k}_1))$ is a disjoint union of
closed stars $\text{Star}(k_s) \, :\, s=1..M$ where $M < N$. In
contrast with the case from the previous paragraph, this time the
restriction of $p$ on each closed disk $\text{Star}({k}_s)$ is a
branch covering map with exactly one ramification point $k_s$ with
index $N_s$ and one branch point $p(k_s) = \tilde{k}_1$.
Consequently, the full preimage of $|\delta|$ is given by the
disjoint union of boundary loops $\partial \text{Star}(k_s)$ for
$s=1..M$. Just like in the previous paragraph, since
$\delta=p(\gamma)$ where $\gamma$ is a connected closed loop in
$\cellcomplex^*(1)$ which bounds an admissible domain on $S$,
$\gamma$ must be among the loops $\partial \text{Star}(k_s)$, i.e.
$\gamma = \partial \text{Star}(k_t)$ for some specific $t \in \{1,
..., M\}$. Consequently the loop $\gamma$ is simple and it splits
the surface $S$ into two open subdomains $\text{OStar}(k_t)$ and
$S \setminus \text{Star}(k_t)$ one of which should be the
admissible domain $\Omega$. Just like before, $S \setminus
\text{Star}(k_t)$ is not a topological disc, but a surface with at
least one handle, so the only option left is $\Omega =
\text{OStar}(k_t)$. As we have assumed that $k_t$ cannot be from
$V_0$ it has to belong to the set of ramification points $V_1$ and
its index of ramification should be $N_{k_t} > 1$. Therefore
$\gamma$, and thus its projection $\delta$, cover the simple loop
$\delta'$ a number of $N_{k_t}$-times. Therefore, by condition 2
of Theorem \ref{thm_Rivin}
$$\sum_{ij^* \in \gamma} (\pi - \theta_{ij})  \, = \,
\sum_{uv^* \in \delta} (\pi - \hat{\theta}_{uv}) \, = \,  N_{k_t}
\sum_{uv^* \in \delta'} (\pi - \hat{\theta}_{uv}) \, = \, 2\pi
N_{k_t} \,
> \, 2\pi.$$

So far we have concluded that whenever $|\delta|$ is simple closed
curve, condition 4 of Schlenker's Theorem \ref{Thm_main} holds.

\medskip

{\bf Subcase 1.2.} Assume $|\delta|$ is not a simple closed curve
on the Riemann sphere, but as a one dimensional connected
subcomplex of $\cellcomplex^*_{\hat{\CC}}(1)$ it is not
simply-connected. Equivalently, $|\delta|$ has a non-trivial
fundamental group. For that reason there exists a cyclic sequence
$\delta'$ of edges of $|\delta|$ that defines a simple closed loop
in the one-skeleton of $\cellcomplex^*_{\hat{\CC}}$. In
particular, $\delta'$ is a cyclic subsequence of $\delta$ and
because $|\delta|$ is not simple while $|\delta'|$ is, $\delta
\setminus \delta' \neq \varnothing$. Consequently, by applying
again condition 2 of Theorem \ref{thm_Rivin}
\begin{align*}
\sum_{ij^* \in \gamma} (\pi-\theta_{ij}) &= \sum_{uv^* \in \delta}
(\pi-\hat{\theta}_{uv})\\ &= \sum_{uv^* \in \delta'}
(\pi-\hat{\theta}_{uv}) + \sum_{uv^* \in \delta \setminus \delta'}
(\pi-\hat{\theta}_{uv}) > \sum_{uv^* \in \delta'}
(\pi-\hat{\theta}_{uv}) \geq 2\pi.
\end{align*}

\medskip

{\bf Subcase 1.3.} $|\delta|$ is a simply-connected subcomplex of
the one-skeleton of the dual complex $\cellcomplex^*_{\hat{\CC}}$.
This is equivalent to saying that $|\delta|$ is a tree in
$\cellcomplex^*_{\hat{\CC}}(1)$ and as such it is contractible to
a point. Therefore, by the lifting properties of covering maps,
the full preimage of $|\delta|$ via $p$ is a disjoint union of $N$
homeomorphic copies of $|\delta|$ in the one-skeleton of the dual
complex $\cellcomplex^*$ on $S$. Since $\gamma$ is a lift of
$\delta$ under $p$ and is connected, $|\gamma|$ should be one of
these copies. Hence, $|\gamma|$ is a tree and furthermore, $\Omega
= S \setminus |\gamma|$. Since $|\gamma|$ is contractible,
$\Omega$ is homeomorphic to $S$ with a closed disc removed, which,
as already pointed out, is not a topological disc. Hence, we
conclude that this case cannot occur either.

To summarize, we have verified condition 4 of Theorem
\ref{Thm_main} whenever the boundary of $\Omega$ consists entirely
of dual edges, i.e. $\partial \Omega \cap V_1 =
\partial \Omega \cap V = \varnothing$.

\paragraph{Case 2.} Let $|\gamma \cap V_1| \geq 2$, which means that
$\gamma$ passes through at least two different vertex points from
$V_1$ or at least twice through the same point from $V_1$. We
claim that since $\gamma$ is the boundary of the admissible domain
$\Omega$, the inequalities $|\gamma \cap V_1| \geq 3$ or $\gamma
\cap E^* \neq \varnothing$ (or both) hold. Indeed, assume that
this is not the case. Then $|\gamma \cap V_1| = 2$ and $\gamma
\cap E^* = \varnothing$. There are only two ways this can happen.
Either $\gamma$ consists of exactly four corner edges or it
consists of two corner edges repeated twice.

In the first case there are exactly two vertices $i$ and $j \in
\gamma \cap \in V_1$, two dual vertices $O_f$ and $O_{f'} \in V^*$
which form the loop of four corner edges $\gamma = (O_fi, iO_{f'},
O_{f'}j, jO_f)$. Both triangles $\hat{\Delta} = iO_fO_{f'}$ and
$\hat{\Delta'} = jO_fO_{f'}$ are two faces of $\hat{\Triang}$ (see
Section \ref{Sec_Existence_Uniqueness_of_Patterns}) that have two
points in common, namely $O_f$ and $O_{f'}$. By strong regularity
of $\hat{\Triang}$ the triangles $\hat{\Delta}$ and
$\hat{\Delta'}$ share a common dual edge $O_fO_{f'}$ and thus
$\gamma$ is the boundary of the topological disc $\hat{\Delta}
\cup \hat{\Delta'}$. Therefore $\gamma$ separates $S$ into two
open subdomains, namely the open interior of $\hat{\Delta} \cup
\hat{\Delta'}$ and $S \setminus (\hat{\Delta} \cup
\hat{\Delta'})$. However, neither of them can be $\Omega$ because
(i) the former is not an admissible domain as it does not contain
any points from $V$ and (ii) the latter is not a topological disc.

In the second case, $\gamma = (O_fi, iO_{f'}, O_{f'}i, iO_f)$.
Then $\Omega$ can only be $S \setminus \gamma$, which is not
possible since $S \setminus \gamma$ is not a topological disc.

Thus, we conclude that $|\gamma \cap V_1| \geq 3$ or $\gamma \cap
E^* \neq \varnothing$, which yields the inequality
$$\sum_{ij^* \subset \gamma} (\pi - \theta_{ij})\,  + \,
\pi|\gamma \cap V_1| \, > \,  2\pi.$$

\paragraph{Case 3.} Let $|\gamma \cap V_1|  = 1$, which means that
$\gamma$ passes through exactly one vertex point $k_1$ from $V_1$
exactly once. Then by projecting down to $\hat{\CC}$ via $p$, we
obtain $\delta = p(\gamma)$ which is a loop in
$\hat{\Triang}_{\hat{\CC}}(1)$ that passes only once through only
one point from $V_1(\hat{\CC})$ denoted by $\tilde{k}_1 = p(k_1)$.
Remove from $\gamma$ the point $k_1$ together with the two corner
edges on $\gamma$ attached to $k_1$ in order to to obtain a path
$\gamma_1$ in the one-skeleton of $\cellcomplex^*$. Let $\delta_1
= p(\gamma_1)$ which is a path in the one-skeleton of
$\cellcomplex^*_{\hat{\CC}}$ obtained by removing $\tilde{k}_1$
and its two adjacent corner edges (which may also be only one
adjacent corner edge repeated twice) from $\delta$.

As already discussed in Subcase 1.1.C, the restriction of $p$ on
the closed star $\text{Star}(k_1)$ is a branch covering map onto
the closed star $\text{Star}(\tilde{k}_1)$ with one ramification
point $k_1$ of ramification index $N_{k_1}$. Let us denote by
$O_{f_1}$ and $O_{f_2}$ the two different dual vertices on the
boundary of $\text{Star}(k_1)$ connected by the path $\gamma_1$.
Then their images $p(O_{f_1})=O_{\tilde{f}_1}$ and
$p(O_{f_2})=O_{\tilde{f}_2}$ are the two dual vertices on the
boundary of $\text{Star}(\tilde{k}_1)$ connected by $\delta_1$.
Just like in the case of loops before, we are going to look at
different cases for the topology of the one dimensional subcomlex
$|\delta_1|$ of $\cellcomplex^*_{\hat{\CC}}(1)$.

\medskip

{\bf Subcase 3.1.} Let $|\delta_1|$ be non-simply connected. This
means that it has a non-trivial fundamental group and so there
exists a cyclic sequence $\delta_1'$ of edges of $|\delta_1|$ that
defines a simple closed loop in $\cellcomplex^*_{\hat{\CC}}(1)$.
Therefore, by condition 2 of Theorem \ref{thm_Rivin}
\begin{align*}
\sum_{ij^* \subset \gamma_1} (\pi - \theta_{ij}) &= \sum_{uv^*
\subset \delta_1} (\pi - \hat{\theta}_{ij}) = \sum_{uv^* \subset
\delta_1'} (\pi - \hat{\theta}_{ij}) + \sum_{uv^*
\subset \delta_1 \setminus \delta_1'} (\pi - \hat{\theta}_{ij}) \\
&\geq 2 \pi + \sum_{uv^* \subset \delta_1 \setminus \delta_1'}
(\pi - \hat{\theta}_{ij}) \, \geq  \, 2 \pi \, > \,  \pi,
\end{align*}
and so $\sum_{ij^* \subset \gamma} (\pi - \theta_{ij}) + \pi > 2
\pi.$

\medskip

{\bf Subcase 3.2.} Assume that $|\delta_1|$ is simply connected.
This means that it is contractible and so it is a tree in the
one-skeleton of $\cellcomplex_{\hat{\CC}}^*$. Then, by the lifting
properties of the regular covering map $p$ on $S \setminus
p^{-1}(V_1(\hat{\CC}))$, the complex $|\gamma_1|$ is also a tree
and is homeomorphic to $|\delta_1|$ via $p$. Consequently, since
the path $\gamma_1$ traverses the tree $|\gamma_1|$, its image
$\delta_1=p(\gamma_1)$ traverses $|\delta_1|$ in the same way.
Adding back the two corner edges $O_{f_1}k_1$ and $k_1O_{f_2}$ to
$\gamma_1$ and $O_{\tilde{f}_1}\tilde{k}_1$ and
$\tilde{k}_1O_{\tilde{f}_2}$ to $\delta_1$ restores the loops
$\gamma$ and $\delta$ respectively, showing that the restricted
map $p|_{|\gamma|} :  |\gamma| \to |\delta|$ is a homeomorphism,
due to the fact that it is a homeomorphism between $|\gamma_1|$
and $|\delta_1|$. Observe that the two corner edges are either
different for each of the two loops $\gamma$ and $\delta$, or they
coincide for each of these two loops.


\medskip

{\bf Subcase 3.2.A.} Assume $O_{\tilde{f}_1} \equiv
O_{\tilde{f}_2}$. This is true exactly when $O_{{f}_1} \equiv
O_{{f}_2}$. Then $|\gamma|$ is contractible, i.e. it is a tree, so
$\Omega = S \setminus |\gamma|$ is not a topological disc when
$g(S) \geq 1$. Hence, this scenario is impossible.

\medskip
{\bf Subcase 3.2.B.} Let $O_{\tilde{f}_1} \neq O_{\tilde{f}_2}$.
This is true exactly when $O_{{f}_1} \neq O_{{f}_2}$. A \emph{true
leaf} of the tree $|\gamma_1|$ is a leaf which is neither
$O_{{f}_1}$ nor $O_{{f}_2}$. A \emph{true leaf-edge} is the unique
edge of the tree attached to a true leaf. The same terminology
applies to $|\delta_1|$.
Remove all true leaves and leaf-edges of $|\gamma_1|$. Perform the
same operation on the homeomorphic tree $|\delta_1|$. Every time
we remove a true leaf together with its corresponding true
leaf-edge from $|\gamma_1|$, we actually add them to the
admissible domain $\Omega$ that $\gamma$ bounds, obtaining a new
admissible domain. As each time we add one vertex with one edge
attached to it, the Euler characteristic of the newly obtained
domain is preserved, i.e. we obtain an open topological disc
again. After removing the true leaves and leaf-edges from
$|\gamma_1|$ and $|\delta_1|$, we end up with a pair of smaller
trees, again homeomorphic via $p$. If these new trees have any
true leaves, we repeat the procedure. We keep repeating until
there are no true leaves left and the only leaves left are
$O_{{f}_1}$ and $O_{{f}_2}$ from $\gamma_1$ and $O_{\tilde{f}_1}$
and $O_{\tilde{f}_2}$ from $\delta_1$. On the level of admissible
domains on $S$, this procedure enlarges the initial admissible
domain $\Omega$ to the admissible domain $\Omega' \, \supset \,
\Omega$, where the latter is obtained by adding to the former all
removed true leaves and leaf-edges. In the end, what is left from
$|\gamma_1|$ and $|\delta_1|$ is a pair of homeomorphic simple
paths in the one-skeletons of the dual complexes $\cellcomplex^*$
and $\cellcomplex^*_{\hat{\CC}}$ respectively. Denote these two
paths by $\gamma_1'$ and $\delta_1'$. Furthermore, $\gamma_1'$
lies on the boundary of $\Omega'$ so that if we add to $\gamma_1'$
the two corner edges $O_{{f}_1}k_1$ and $k_1O_{{f}_2}$ then we
obtain the full boundary, call it $\gamma'$, of $\Omega'$.
Consequently, $\gamma'$ is a simple loop in the one-skeleton of
$\cellcomplex^*$ and $\gamma' = \partial \Omega'$. Recall that by
construction $\delta'_1=p(\gamma'_1)$ is a simple path and so
$\delta'=p(\gamma')$ is a simple loop.

 Assume that $\delta'_1$ lies on the boundary of a closed star
 $\text{Star}(\tilde{k})$ for some $\tilde{k} \in V_{\hat{\CC}}$.
 Since, by construction, the cell complexes on $S$ are the
 lifts of the cell complexes on $\hat{\CC}$ via $p$, the
 preimage of $\text{Star}(\tilde{k})$ is a
 disjoint union of stars and thus $\gamma_1'$ lies in the boundary
 of a star $\text{Star}({k})$ for some $k \in V$ such that
 $p(k) = \tilde{k}$.

 First, if $k \neq k_1$, then the endpoints $O_{f_1}$ and $O_{f_2}$ of
 $\gamma_1'$ lie simultaneously in $\text{Star}(k)$ and
 $\text{Star}(k_1)$. As these two endpoints are by assumption
 different, the two stars share two different vertices
 and so, by strong regularity of $\hat{\Triang}$, they share
 exactly one common dual edge $O_{f_1}O_{f_2} \in E^*$. Then, if
 we denote by $\hat{\Delta}$ the (closed) triangular face
 $k_1O_{f_1}O_{f_2}$ of $\hat{\Triang}$, we conclude that the simple loop
 $\gamma'$ splits the surface $S$ into two open domains, one of
 which is the open interior of $\text{Star}(k) \cup \hat{\Delta}$
 and the other is $S \setminus \big(\text{Star}(k) \cup
 \hat{\Delta}\big)$. None of them can be $\Omega'$ because the
 latter is not an open topological disc, whenever $g(S) \geq 1$,
 while the former is a topological disc, but it is not an admissible domain,
 as it cannot be represented as the union of open stars of
 $\hat{\Triang}$ due to the presence of the
 additional triangle $\hat{\Delta}$. Therefore, this situation
 cannot occur.

Second, if $k = k_1$, then on one side $\gamma'$ bounds on $S$ a
strict subdomain of $\text{Star}(k_1)$ which does not even contain
an open star, so the domain cannot be admissible, while on the
other side of $\gamma'$ we have $S$ with a closed topological disc
removed, which, as before, is not an open topological disc.
Therefore, this situation cannot occur either.

Consequently, $\delta'_1$ cannot lie on the boundary of any star,
which is equivalent to saying that the simple path of dual edges
$\delta'_1$ does not lie in the boundary of a dual face of
$\cellcomplex^*_{\hat{\CC}}$. By Bao and Bonahon's Theorem
\ref{thm_Bao} $\sum_{ij^* \subset \gamma_1'} (\pi - \theta_{ij}) =
\sum_{uv^* \subset \delta_1'} (\pi - \hat{\theta}_{uv}) > \pi.$
Returning to the original loop $\gamma = \partial \Omega$ which
contains $\gamma'_1$
\begin{align*}
\sum_{ij^* \subset \gamma} (\pi - \theta_{ij}) + \pi &= \sum_{ij^*
\subset \gamma_1'} (\pi - \theta_{ij}) + \sum_{ij^* \subset \gamma
\setminus \gamma_1'} (\pi - \theta_{ij}) + \pi \\
&> \pi + \sum_{ij^* \subset \gamma \setminus \gamma_1'} (\pi -
\theta_{ij}) + \pi \geq 2 \pi.
\end{align*}

Finally, we have concluded that condition 4 of Theorem
\ref{Thm_main} holds for all possible cases. Thus, the proof of
Theorem \ref{thm_general_uniformization} is complete.


\section{Proof of Theorem \ref{thm_sphere_uniformization}}
\label{Sec_proof_sphere_uniformization}

The proof of Theorem \ref{thm_sphere_uniformization} is analogous
to the proof of Theorem \ref{thm_general_uniformization} from the
preceding section. The major difference is that instead of
Schlenker's Theorem \ref{Thm_main}, our main tool is Bao and
Bonahon's Theorem \ref{thm_Bao} because $S$ is a topological
sphere. Consequently, instead of working with strict admissible
domains with different topologies, we check the necessary and
sufficient conditions of Theorem \ref{thm_Bao} for simple loops
$\gamma$ in the one-skeleton of the dual complex $\cellcomplex^*$
on $S$, as well as for simple paths $\gamma_1$ in
$\cellcomplex^*(1)$ joining two dual vertices from the same dual
face without being contained in the boundary of any dual face. In
the case of loops, we simply repeat the arguments from Subcases
1.1. and 1.2 in Section
\ref{Sec_proof_uniformization_branch_covers}, observing that due
to the simplicity of $\gamma$, Subcase 1.3 cannot occur. In the
case of paths, we repeat the arguments from Subcases 3.1 and 3.2
with the simplification that neither $|\gamma_1|$ nor
$p(|\gamma_1|)$ have true leaves, so the procedure of removing
them is unnecessary in this case.




\section{Spaces of decorated triangles and circle
patterns} \label{Sec_Spaces_of_triangles_and_patterns}

In order to outline the algorithmic recipe for discrete
uniformzation we need to provide some background on the spaces of
decorated triangles and hyper-ideal circle patterns. In sections
\ref{Sec_uniformization_discrete_negative_surfaces} and
\ref{Sec_uniformization_hyper_elliptic_surfaces} we have explained
how one can obtain in a natural and systematic way combinatorial
angle data $(\cellcomplex, \theta, \Theta)$ from either a
negatively curved or a hyper-elliptic Riemann surface.
Theorem \ref{Thm_main} guarantees the unique geometric
realizability of the data $(\cellcomplex, \theta, \Theta)$. In the
following two sections \ref{Sec_Spaces_of_triangles_and_patterns}
and \ref{Sec_Variational principle} we describe the variational
principle which plays a central role in the verification of
Theorem \ref{Thm_main} (see \cite{ND}) and provides a method for
the construction of hyper-ideal circle patterns.
All constructions and notations follow closely the exposition of
\cite{ND}. The interested reader can check the necessary details
there. In the current article, we provide a sketch.

\begin{figure}
\centering
\includegraphics[width=9.5cm]{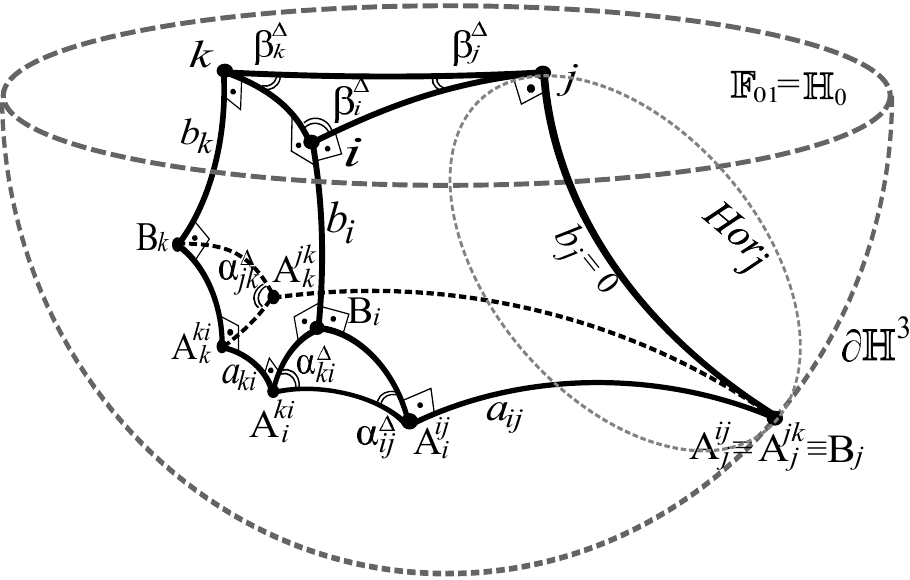}
\caption{A hyper-ideal tetrahedron $\tau_{\Delta}$.} \label{Fig6}
\end{figure}

\subsection{The space of decorated triangles}

First, we describe the space of (labelled) decorated triangles in
$\hyperbolicplane$ considered up to hyperbolic isometries, with
predetermined fixed splitting of the triangles' vertices
$V_{\Delta} = \{i,j,k\}$ into $V_{\Delta}^1$ and $V_{\Delta}^0$.
Then, up to $\hyperbolicplane$-isometry, a decorated triangle
$\Delta=ijk$ can be uniquely represented in three different ways
\cite{ND}. The most natural way is to give the triangles edge
lengths and vertex radii $(l,r)_{\Delta} = (l_{ij}, l_{jk},
l_{ki}, r_k, r_i, r_j) \in \ERD$, where the first three are
positive and satisfy all three strict triangle inequalities
$l_{uv} < l_{vw} + l_{wu}$, as well as $l_{uv} > r_u + r_v$ for
all $u \neq v \neq w \in \{i,j,k\}$, while the last three satisfy
$r_u > 0$ for $u \in V_{\Delta}^1$ and $r_u=0$ for $u \in
V_{\Delta}^0$. The second way of uniquely representing a decorated
triangle is by its six angles $(\alpha^{\Delta}, \beta^{\Delta}) =
(\alpha_{ij}^{\Delta}, \alpha_{jk}^{\Delta}, \alpha_{ki}^{\Delta},
\beta_{k}^{\Delta}, \beta_{i}^{\Delta}, \beta_{j}^{\Delta}) \in
\ADelta$ \cite{Sch1,S,ND}. Here $\alpha^{\Delta}_{uv}$ is the
angle between the geodesic edge $uv$ of $\Delta=ijk$ and its face
circle $c_{\Delta}$ measured inside $c_{\Delta}$ and outside
$\Delta$, and $\beta^{\Delta}_v$ is the interior angle of the
triangle at its vertex $v \in V_{\Delta}$ (refer to Figure
\ref{Fig1}). The third type of decorated triangle description is
the one that is most crucial to our later constructions. It draws
upon the link between decorated triangles in $\hyperbolicplane$
and hyper-ideal tetrahedra in $\hyperbolicspace$.

\subsection{Hyper-ideal tetrahedra}

\begin{Def} \label{Def_hyperideal_tetrahedron}
A \emph{hyper-ideal tetrahedron} (see \cite{Sch1,S} and Figure
\ref{Fig6}) is a geodesic polyhedron in $\hyperbolicspace$ that
has the combinatorics of a tetrahedron with some (possibly all) of
its vertices truncated by triangular \emph{truncating faces}. Each
truncating face is orthogonal to the faces and the edges it
truncates. Furthermore, a pair of truncating faces do not
intersect. Finally, the non-truncated vertices are all ideal.
\end{Def}

The polyhedron depicted on Figure \ref{Fig6} is an example of a
hyper-ideal tetrahedron with three truncated \emph{hyper-ideal
vertices} and one ideal vertex. This terminology comes from the
interpretation that in the Klein projective model or the Minkowski
space-time model of $\hyperbolicspace$ \cite{ThuBook,BenPetr,S}
the hyper-ideal tetrahedron can be represented by an actual
tetrahedron with some (or all) vertices lying outside
$\hyperbolicspace$ (hence the term \emph{hyper-ideal vertices}).
The dual (projective polar) to each hyper-ideal vertex is the
orthogonal truncating plane. However, in this article, we mostly
use the two standard conformal models of $\hyperbolicspace$ - the
Poincar\'e ball model and the upper half-space model
\cite{ThuNotes,ThuBook,BenPetr}, which are the three dimensional
analogs of the disk model and the upper half-plane model of
$\hyperbolicplane$ respectively. 
We call \emph{a principal edge} a geodesic edge of a hyper-ideal
tetrahedron that is not contained in a truncating plane of the
hyper-ideal tetrahedron. Each hyper-ideal tetrahedron has exactly
six of them. The rest of the edges lie on the truncating planes.
We can call them \emph{auxiliary edges}. A principal edge has one
of the following three properties: (i) it goes from one truncating
face to another, being by definition orthogonal to both of them
(think of it as the two-sided truncation of the edge connecting
two hyper-ideal vertices); (ii) it goes from one ideal vertex to a
truncating face, being by definition orthogonal to the truncating
face (think of it as the one-sided truncation of the edge
connecting an ideal vertex to a hyper-ideal one); (iii) it
connects two ideal vertices. See Figure \ref{Fig6}.

\begin{figure} \centering
\includegraphics[width=13cm]{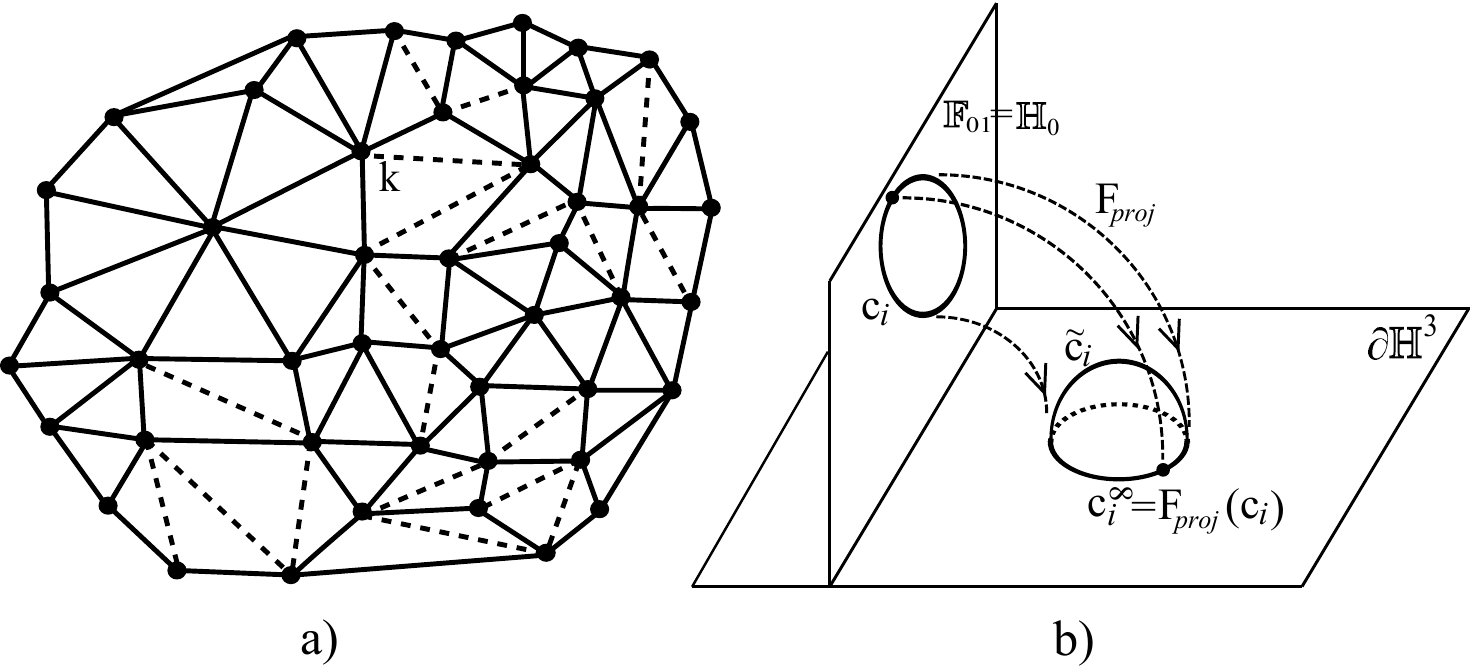}
\caption{a) The subtriangulation $\Triang=(V,E_T,F_T)$ of
$\cellcomplex$ whose dashed edges are the auxiliary edges from
$E_{\pi}$; b) The projection $F_{proj}$ from the hyperbolic plane
$\Hplane \subset \hyperbolicspace$ to the ideal boundary
$\partial\hyperbolicspace$.} \label{Fig5}
\end{figure}

\subsection{Link between decorated triangles and hyper-ideal
tetrahedra}

To obtain a hyper-ideal tetrahedron $\tau_{\Delta}$ from a
decorated triangle $\Delta=ijk$, first think of
$\hyperbolicplane$, together with the decorated triangle $\Delta$
drawn on it, as a hyperbolic plane $\Hplane \cong
\hyperbolicplane$ lying in $\hyperbolicspace$. This situation is
depicted on Figure \ref{Fig5}. The idea is that we can project
$\Hplane$ down on $\partial\hyperbolicspace$ using a natural
geometric map $F_{proj} \, : \, \Hplane \, \to \,
\partial\hyperbolicspace$ (see also Figure \ref{Fig5}). For any
$x \in \Hplane$ take the unique geodesic in $\hyperbolicspace$
passing through $x$ and perpendicular to $\Hplane$ and follow it,
only on one side of $\Hplane$, all the way down to
$\partial\hyperbolicspace$ reaching the ideal point $F_{proj}(x)$.
Then the image $F_{proj}(\Delta)$ of the decorated triangle
$\Delta$ determines a set of circles (and straight lines,
depending on the model) in $\partial\hyperbolicspace$ which when
extended to hyperbolic planes cut out a hyper-ideal tetrahedron.
Figure \ref{Fig5} shows how a vertex circle $c_i$ is being mapped
to a circle $F_{proj}(c_i) = c_i^{\infty}$ on
$\partial\hyperbolicspace$ and then extended to a hyperbolic plane
$\tilde{c}_i$ of $\hyperbolicspace$. The converse construction
also holds in the sense that a (labelled) hyper-ideal tetrahedron
$\tau$ gives rise to a unique decorated triangle
$\triangle_{\tau}$. The face labelled $ijk$ determines the
hyperbolic plane $\Hplane$. Then, we can take all the ideal
circles (and possibly straight lines, depending on the
$\hyperbolicspace$ model) on $\partial\hyperbolicspace$ of all
hyperbolic planes that the faces of $\tau$ determine, and map them
back to $\Hplane$ via $F_{proj}^{-1}$. Notice that since all
constructions utilize only the geometry of $\hyperbolicspace$,
they are invariant with respect to hyperbolic congruences (that
respect the labelling of the triangles and the tetrahedra).
Furthermore, because the models of both $\hyperbolicplane \cong
\Hplane$ and $\hyperbolicspace$ are conformal, the six angles
$(\alpha^{\Delta},\beta^{\Delta})$ of the decorated triangle
become the corresponding dihedral angles at the six principal
edges of the constructed hyper-ideal tetrahedron $\tau_{\Delta}$.
This can be seen on Figure \ref{Fig6}.

\subsection{Tetrahedral edge-length variables}

Now, after we have explained the one-to-one correspondence between
(labelled) decorated triangles in $\hyperbolicplane$ and
(labelled) hyper-ideal tetrahedra, we can define the variables
$(a,b)_{\Delta} = (a_{ij}, a_{jk}, a_{ki}, b_k, b_i, b_j) \in
\TED$ as the ``hyperbolic lengths" of the principal geodesic edges
of $\tau_{\Delta}$. For principal edges of type (i), as described
above (i.e. edges perpendicular to two truncating faces), the
geodesic length makes sense and is a positive real number. But
edges of type (ii) and (iii) have actually infinite hyperbolic
length. In order to fix this issue, we \emph{decorate} the
hyper-ideal tetrahedron $\tau_{\Delta}$, which was obtained from
the decorated triangle $\Delta=ijk$, with one horosphere per ideal
vertex so that the horosphere touches $\partial\hyperbolicspace$
at that ideal vertex and the plane $\Hplane$ at the endpoint of
the geodesic edge emanating from the ideal vertex in question. On
Figure \ref{Fig6} $B_j$ is an ideal vertex and the decorating
horosphere $\mathit{Hor}_j$ is tangent to
$\partial\hyperbolicspace$ at $B_j$ and to $\Hplane$ at $j$.
Consequently, the length of a type (ii) geodesic edge of
$\tau_{\Delta}$ is the oriented hyperbolic distance between the
truncating face on one side of the edge and the decorating
horosphere on the other side, measured along the edge itself. The
length is positive if the truncating face and the decorating
horosphere are disjoint, zero if they are tangent and negative if
they intersect. Similarly, the length of a type (iii) geodesic
edge of $\tau_{\Delta}$ is the oriented hyperbolic distance
between the two decorating horospheres, one on each side of the
edge, measured along the edge itself. As before, the length is
positive if the two decorating horospheres are disjoint, zero if
they are tangent and negative if they intersect. Edges of type (i)
and (ii) can be seen on Figure \ref{Fig6}. For instance, $a_{ki},
b_k$ and $b_i$ are the hyperbolic lengths of geodesic edges of
type (i), while $a_{ij}$ is the length of the edge $A^{ij}_iB_j$
defined as the (oriented) distance between the truncating face
$A^{ki}_iA^{ij}_iB_i$ and the decorating horosphere
$\mathit{Hor}_j$. Notice that $b_j=0$ since $\text{Hor}_j$ is by
construction tangent to $\Hplane$ at the point $j$.

\subsection{Transition formulas between different sets of
variables}

It is very important to find how the angles
$(\alpha^{\Delta},\beta^{\Delta}) \in \ADelta$ of a decorated
triangle $\Delta$, which are also the six principal dihedral
angles of the corresponding hyper-ideal tetrahedron
$\tau_{\Delta}$, depend on the principal edge-lengths
$(a,b)_{\Delta} \in \TED$ of $\tau_{\Delta}$. It is also useful to
know how the edge-lengths and vertex radii $(l,r)_{\Delta} \in
\ERD$ of $\Delta$ depend on the parameters $(a,b)_{\Delta}$. By
using various combinations of hyperbolic trigonometric formulas
\cite{Bus} applied to the faces of $\tau_{\Delta}$ one can derive
the necessary expressions. In what follows, we present in detail
the formulas for hyperbolic decorated triangles, as this is the
more general case of higher genus surfaces. The Euclidean case is
analogous and can be worked out by applying the appropriate
formulas given in \cite{ND} for instance. For the edge-lengths and
vertex radii of a hyperbolic $\Delta=ijk$ the following formulas
apply \cite{ND}:
\begin{align}
r_v &= \sinh^{-1}{\left(\frac{1}{\sinh{b_v}}\right)} \,\text{ if }
\, v \in V_{\Delta}^1 \,\,\,\, \text{ and } \,\,\,\, r_v=b_v=0
\,\, \text{ if } \, v
\in V^0_{\Delta} \label{Eqn_Hyp_r_b}\\
l_{uv} &= \cosh^{-1}{\left(\frac{\cosh{b_u} \cosh{b_v}
+\cosh{a_{uv}}}{\sinh{b_u}\sinh{b_v}}\right)} =
f_3(b_u,b_v,a_{uv})
\,\, \text{ if } \, u,v \in V^1_{\Delta} \label{Eqn_Hyp_l_a_b_b}\\
l_{uv} &= \cosh^{-1}{\left(\frac{\cosh{b_u} +
e^{a_{uv}}}{\sinh{b_u}}\right)} = f_2(b_u,a_{uv}) \,\,\, \text{ if
}
\, u \in V^{\Delta}_1, \,\text{ and } \, v \in V^{\Delta}_0  \label{Eqn_Hyp_l_a_b}\\
l_{uv} &= 2 \sinh^{-1}{\big(e^{a_{uv}/2}\big)} = f_1(a_{uv})
\,\,\, \text{ if } \, u,v \in V^0_{\Delta}, \label{Eqn_Hyp_l_a}
\end{align}
where $uv \in E_{\Delta} = \{ij, jk, ki\}$ is an edge of
$\Delta=ijk$. Having computed the edge-lengths of $\Delta$ one can
immediately find, using for example the hyperbolic law of cosines,
the angles
\begin{align} \label{Eqn_angles_beta}
\beta_v^{\Delta} = \arccos\left(\frac{\cosh{l_{uv}}\cosh{l_{vw}} -
\cosh{l_{wu}}}{\sinh{l_{uv}}\sinh{l_{vw}}}\right)=g(l_{uv},l_{vw},l_{wu})
\,\, \text{ for all } \, v \in V_{\Delta}.
\end{align}
In order to find the angles $\alpha_{uv}^{\Delta}$ we examine all
four combinatorial types of decorated triangle. There are
different ways and different hyperbolic trigonometric formulas one
can use, but we show only some of them (all leading to the same
result).

\paragraph{Case 1.} $V_{\Delta}^1 = \{i,k\}$ and $V_{\Delta}^0=\{j\}$.
The hyper-ideal tetrahedron $\tau_{\Delta}$ corresponding to this
combinatorics is depicted on Figure \ref{Fig6}. We work with the
notations introduced there for the current case as well as for the
rest of the cases. Compute the geodesic edge-lengths of the
triangular truncating face $A^{ki}_iA^{ij}_iB_i$
\begin{align*}
\sigma^{ki}_i &= l_{\hyperbolicspace}\big(A^{ki}_iB_i\big) =
\cosh^{-1}{\left(\frac{\cosh{a_{ki}} \cosh{b_i}
+\cosh{b_k}}{\sinh{a_{ki}}\sinh{b_i}}\right)} =
f_3(a_{ki},b_i,b_k)\\
\sigma^{ij}_i &= l_{\hyperbolicspace}\big(A^{ij}_iB_i\big) =
\cosh^{-1}{\left(\frac{\cosh{b_i} +
e^{-a_{ij}}}{\sinh{b_i}}\right)} = f_2(b_i,-a_{ij}) \\
\sigma_i &= l_{\hyperbolicspace}\big(A^{ki}_iA^{ij}_i\big) =
\cosh^{-1}{\left(\frac{\cosh{a_{ki}} +
e^{a_{jk}-a_{ij}}}{\sinh{a_{ki}}}\right)} =
f_2(a_{ki},a_{jk}-a_{ij}).
\end{align*}
Then $\, \alpha_{ij}^{\Delta} = g\big(\sigma_i, \sigma^{ij}_i,
\sigma^{ki}_i\big) \,$ and $ \, \alpha_{ki}^{\Delta} =
g\big(\sigma_i, \sigma^{ki}_i, \sigma^{ij}_i\big) \,$ (see formula
(\ref{Eqn_angles_beta})). Alternatively, one can also use the
hyperbolic law of sines
\begin{align*}
\alpha_{ij}^{\Delta} =
\arcsin\left(\frac{\sinh{\sigma^{ki}_i}}{\sinh{\sigma_i}}
\sin{\beta^{\Delta}_i}\right) \,\,\, \text{ and } \,\,\,
\alpha_{ki}^{\Delta} =
\arcsin\left(\frac{\sinh{\sigma^{ij}_i}}{\sinh{\sigma_i}}
\sin{\beta^{\Delta}_i}\right)
\end{align*}
The easiest way to compute the last angle is $\alpha^{\Delta}_{jk}
= \pi - \alpha^{\Delta}_{ij} - \beta^{\Delta}_{j}$.

\paragraph{Case 2.} $V_{\Delta}^1 = V_{\Delta}=\{i,j,k\}$ while
$V_{\Delta}^0=\varnothing$. One way of computing the angles is to
apply the hyperbolic law of cosines to two out of the three
truncating faces different from the truncating face $ijk$. For
instance, compute the geodesic edge-lengths of the triangular
truncating faces $A^{ki}_iA^{ij}_iB_i$ and $A^{ij}_jA^{jk}_jB_j$.
Obtain
\begin{align*}
\sigma^{ij}_i &= f_3\big(a_{ij}, b_i, b_j\big), &\sigma^{ki}_i &=
f_3\big(a_{ki}, b_i, b_k\big), &\sigma_i &= f_3\big(a_{ki},
a_{ij},
a_{jk}\big) &\text{for }& A^{ki}_iA^{ij}_iB_i\\
\sigma^{ij}_j &= f_3\big(a_{ij}, b_j, b_i\big), &\sigma^{jk}_j &=
f_3\big(a_{jk}, b_j, b_k\big), &\sigma_j &= f_3\big(a_{ij},
a_{jk}, a_{ki}\big)  &\text{for }& A^{ij}_jA^{jk}_jB_j.
\end{align*}
Then $\, \alpha_{ij}^{\Delta} = g\big(\sigma_i, \sigma^{ij}_i,
\sigma^{ki}_i\big) = g\big(\sigma_j, \sigma^{ij}_j,
\sigma^{jk}_j\big) \,$, as well as $ \, \alpha_{jk}^{\Delta} =
g\big(\sigma_j, \sigma^{jk}_j, \sigma^{ij}_j\big) \,$  and $ \,
\alpha_{ki}^{\Delta} = g\big(\sigma_i, \sigma^{ki}_i,
\sigma^{ij}_i\big) \,$ (see formula (\ref{Eqn_angles_beta})).
Alternatively, just like in Case 1, one can also use the
hyperbolic law of sines, applied to $A^{ki}_iA^{ij}_iB_i$ and
$A^{ij}_jA^{jk}_jB_j$ whose angles $\beta^{\Delta}_i$ and
$\beta^{\Delta}_j$ we already know.

\paragraph{Case 3.} $V_{\Delta}^1 = \{i\}$ and $V_{\Delta}^0=\{j, k\}$.
In this case, for example, compute the edge-lengths of the
truncating face $A^{ki}_iA^{ij}_iB_i$
\begin{align*}
\sigma^{ij}_i &= f_2\big(b_i, -a_{ij}\big), &\sigma^{ki}_i &=
f_2\big(b_i, -a_{ki}\big), &\sigma_i &= f_1\big(a_{ki} + a_{ij} +
a_{jk}\big) & \,\,\,\, &
\end{align*}
and then calculate $\, \alpha_{ij}^{\Delta} = g\big(\sigma_i,
\sigma^{ij}_i, \sigma^{ki}_i\big)$. After that
$\alpha^{\Delta}_{jk} = \pi - \beta_j - \alpha^{\Delta}_{ij}$ and
$\alpha^{\Delta}_{ki} = \alpha^{\Delta}_{ij} + \beta_j -
\beta^{\Delta}_{k} = g\big(\sigma_i, \sigma^{ki}_i,
\sigma^{ij}_i\big)$.

\paragraph{Case 4.} $V_{\Delta}^1 = \varnothing$ and
$V_{\Delta}^0=V_{\Delta} = \{i, j, k\}$. Then
\begin{align*}
\alpha^{\Delta}_{ij} &= \frac{\pi + \beta^{\Delta}_k -
\beta^{\Delta}_i - \beta^{\Delta}_j}{2}, &\alpha^{\Delta}_{jk} &=
\frac{\pi + \beta^{\Delta}_i - \beta^{\Delta}_j -
\beta^{\Delta}_k}{2}, &\alpha^{\Delta}_{ki} &= \frac{\pi +
\beta^{\Delta}_j - \beta^{\Delta}_k - \beta^{\Delta}_i}{2}.
\end{align*}


\subsection{The space of generalized hyper-ideal circle
patterns}

As already mentioned, we construct the hyper-ideal circle pattern
that realizes the data $(\cellcomplex,\theta,\Theta)$ as the
unique critical point of a convex functional defined on a suitable
space of patterns. In what follows, we define this space.

\paragraph{Description in terms of edge-lengths and vertex radii of decorated triangles.}
Let $\Triang = (V, E_T, F_T)$ be a subtriangulation of
$\cellcomplex$ obtained by adding a maximal number of diagonals
with non-intersecting interiors in each non-triangular face of
$\cellcomplex$ (see Figure \ref{Fig5}). Consequently, the elements
of $F_T$  are all combinatorial triangles and $E_T = E \cup
E_{\pi}$, where $E_{\pi}$ is the set of all diagonals we have
added in the process of subtriangulation of $\cellcomplex$. On
Figure \ref{Fig5} these are all the dashed edges, while all solid
edges are the elements of $E$. We refer to the edges from
$E_{\pi}$ as \emph{redundant edges}. With this new combinatorics
at hand, define the space of all \emph{generalized} hyper-ideal
circle patterns with combinatorics $\Triang$ on the surface $S$,
considered up to isometry isotopic to identity as follows: the
vector $(l,r) \in \reals^{E_T} \times \reals^V$ belongs to $\ER$
if and only if

\smallskip

\noindent $\bullet$ $l_{ij} > 0$ for $ij \, \in \, E_T$, as well
as $\, r_k > 0$ for $k \in V_1$ and $\, r_k = 0$ for $k \in V_0
\,$;

\smallskip

\noindent $\bullet$ $l_{ij} > r_i + r_j$ for $ij \in E_T$;

\smallskip

\noindent $\bullet$ $l_{ij} < l_{jk} + l_{ki}, \,\,\,\,  l_{jk} <
l_{ki} + l_{ij}, \,\,\,\, l_{ki} < l_{ij} + l_{jk}$ for
$\Delta=ijk \in F_T$.

\smallskip

By definition, the space $\ER$ is clearly a convex polytope in
$\reals^{E_T} \times \reals^{V}$ of dimension $\dim \ER = |E_T| +
|V_1|.$ We can consider the vector space $\reals^{E_T} \times
\reals^{V_1}$ as the vector subspace of $\reals^{E_T} \times
\reals^{V}$ defined by setting $r_k=0$ for all $k \in V_0$. Then
it is clear that $\ER$ is an open convex polytope of $\reals^{E_T}
\times \reals^{V_1}$. The term \emph{generalized hyper-ideal
circle pattern} is used because the hyper-ideal circle patterns
that satisfy the conditions of $\ER$ do not necessarily satisfy
the local Delaunay property from Definition
\ref{Def_local_Del_property}.

\paragraph{Description in terms of edge-lengths of hyper-ideal tetrahedra.} Furthermore, one can define the set $\TE \subset \reals^{E_T}
\times \reals^{V},$ which is in fact an open subset of
$\reals^{E_T} \times \reals^{V_1},$ as the domain of a
real-analytic map $\Psi \, : \, \TE \, \to \, \ER$ defined by
formulas (\ref{Eqn_Hyp_r_b}) to (\ref{Eqn_Hyp_l_a}) so that
$\Psi(a,b) = (l,r)$ and $\Psi(\TE)=\ER$. It is straightforward to
verify that formulas (\ref{Eqn_Hyp_r_b}) to (\ref{Eqn_Hyp_l_a})
can be inverted and an inverse map $\Psi^{-1}\, : \, \ER \, \to \,
\TE$ can be obtained, which is also real-analytic. Thus, one sees
that $\Psi$ is a real-analytic diffeomorphism between the open
subsets $\TE$ and $\ER$ of $\reals^{E_T} \times \reals^{V_1}$.
Therefore, $\TE$ also defines the space of all generalized
hyper-ideal circle patterns with combinatorics $\Triang$ on $S$,
considered up to hyperbolic isometries isotopic to identity. 

The hyper-ideal circle pattern that realizes the data
$(\cellcomplex, \theta, \Theta)$ can be subtriangulated by adding
all geodesic redundant edges from $E_{\pi}$ so that now it has
combinatorics $\Triang = (V, E_T, F_T)$ instead of $\cellcomplex =
(V, E, F)$, where recall that $E_T = E \cup E_{\pi}$.
Consequently, $ij \in E_{\pi}$ if and only if the two compatibly
adjacent decorated triangles that share $ij$ as a common geodesic
edge have coinciding face circles, i.e. they share the same face
circle. This is equivalent to the fact that the intersection angle
between the two adjacent face circles is $\theta_{ij} = \pi$.
Consequently, one can extend $\theta \, : \, E \, \to \, (0, \pi)$
to $\tilde\theta \, : \, E_T \, \to \, (0, \pi]$ by setting
$\tilde{\theta}_{ij}= \theta_{ij}$ whenever $ij \in E$ and
$\tilde{\theta}_{ij}= \pi$ whenever $ij \in E_{\pi}$. Then the
pattern which realizes $(\Triang, \tilde{\theta}, \Theta)$ is
exactly the pattern which realizes the original data
$(\cellcomplex, \theta, \Theta)$ after erasing the redundant edges
$E_{\pi}$. Thus the pattern we are looking for can be seen as a
special pattern lying inside the space $\TE$.

\section{Variational principle for construction of circle
patterns} \label{Sec_Variational principle}

In this section we construct the functional whose only critical
point is the unique pattern that realizes the combinatorial angle
data $(\Triang, \tilde{\theta}, \Theta)$, which as discussed in
the previous Section \ref{Sec_Spaces_of_triangles_and_patterns},
is also the pattern that realizes $(\cellcomplex, \theta, \Theta)$
after the removal of the redundant edges.



Take any $\Delta \in F_T$. For each $(a,b)_{\Delta} \in \TED$
define the function
\begin{align} \label{Formula_local_functional}
\UDelta(a,b) &= \sum_{ij \in E_{\Delta}} \alpha^{\Delta}_{ij}
a_{ij} + \sum_{k \in V_{\Delta}} \beta^{\Delta}_{ij} b_{k} + 2
\Vol(\alpha^{\Delta}, \beta^{\Delta}) \nonumber \\
&= \sum_{ij \in E_{\Delta}} \alpha^{\Delta}_{ij} a_{ij} + \sum_{k
\in V_{\Delta}^1} \beta^{\Delta}_{ij} b_{k} + 2
\Vol(\alpha^{\Delta}, \beta^{\Delta})
\end{align}
where, as already discussed in Section
\ref{Sec_Spaces_of_triangles_and_patterns} (cases 1 to 4), the
angles $\alpha^{\Delta}_{ij}=\alpha^{\Delta}_{ij}(a,b)$ for $ij
\in E_{\Delta}$ and $\beta^{\Delta}_k=\beta^{\Delta}_k(a,b)$ for
$k \in V_{\Delta}$ are real-analytic functions depending on the
tetrahedral edge-length variables $(a,b)_{\Delta} \in \TED$.
Recall that $E_{\Delta}$ is the set of edges of $\Delta$ and
$V_{\Delta}$ is the set of its vertices with $V_{\Delta}^1$ being
the subset of those vertices of $\Delta$ that are supposed to have
vertex circles of positive radius. The vertices from its
complement $V_{\Delta}^0$ satisfy the restriction $b_k = 0, \, k
\in V_{\Delta}^0$, which leads us to the second sum in
(\ref{Formula_local_functional}). The function $\Vol$, which
depends analytically on the angles $(\alpha^{\Delta},
\beta^{\Delta})$, is the hyperbolic volume of the hyper-ideal
tetrahedron $\tau_{\Delta}$ with principal edge-lengths
$(a,b)_{\Delta} \in \TED$ and corresponding dihedral angles
$(\alpha^{\Delta}, \beta^{\Delta})$.

As a function of the dihedral angles, $\Vol$ is strictly concave
\cite{Sch1,Sch3,S} and because of that, as shown in \cite{ND},
$\UDelta$ is a locally strictly convex function on $\TED$. It is
straightforward to verify that each real-analytic angle function
$\alpha_{ij}^{\Delta} = \alpha^{\Delta}_{ij}(a,b)$ and
$\beta_{k}^{\Delta} = \beta^{\Delta}_{k}(a,b)$ can be continuously
extended by $\alpha^{\Delta}_{ij}\equiv \pi \equiv
\beta^{\Delta}_k$ whenever $l_{ij} \geq l_{jk} + l_{ki}$, and
$\alpha^{\Delta}_{ij}\equiv 0 \equiv \beta^{\Delta}_k$ whenever
either $l_{jk} \geq l_{ki} + l_{ij}$ or $l_{ki} \geq l_{ij} +
l_{jk}$. Let us partition the set $E_{\Delta}$ of edges of
$\Delta$ into
$$E_{\Delta}^{1} = \{ij \in
E_{\Delta} \, | \, i, j \in V_{\Delta}^1\}  \, \text{ and } \,
E_{\Delta}^{0} = E_{\Delta} \, \setminus \, E_{\Delta}^1.$$ Then
the angle functions $\alpha^{\Delta}$ and $\beta^{\Delta}$ are
continuous on $\reals^{E_{\Delta}^0} \times
\reals_{+}^{E_{\Delta}^1 \cup V_{\Delta}^1}$. Furthermore, outside
$\TED$, the volume $\Vol$ is constantly zero. Consequently, as
explained for example in \cite{BPS}, the function $\UDelta$, which
is real-analytic and locally strictly convex in $\TED$, can be
extended to a continuously differentiable convex function on the
whole convex set $\reals^{E_{\Delta}^0} \times
\reals_{+}^{E_{\Delta}^1 \cup V_{\Delta}^1} \, \supset \TED$.

The reason for which $\UDelta$ is suitable for applications is
that there exist formulas for the volume $\Vol$ of a hyper-ideal
tetrahedron with given principal dihedral angles
$(\alpha^{\Delta},\beta^{\Delta})$. For instance, due to
Springborn \cite{S}, in the case of a hyper-ideal tetrahedron with
at least one ideal vertex there is a fairly nice explicit
expression for its volume in terms of Lobachevsky's functions.
This formula however does not work for a hyper-ideal tetrahedron
with exactly four hyper-ideal vertices. Nevertheless, a formula in
that case also exists due to Ushijima \cite{Ush}. Moreover, it
cover all possible combinatorial types of hyper-ideal tetrahedra.
However, it is quite more complicated and in our computer
realizations of the discrete uniformization procedure, whose
results can be found in Section \ref{Sec_computer_realizations},
we have mostly used Springborn's formula (and its
simplifications), while Ushijima's version have been used only for
the cases of tetrahedra with four hyper-ideal vertices. All of
these formulas are real-analytic in nature.


With the data $(\Triang, \tilde{\theta}, \Theta)$ at hand, one can
construct the functional
\begin{align} \label{Formula_global_functional_general}
\Uglobal_{\theta, \Theta}(a,b) = \sum_{\Delta \in F_T}
\UDelta(a,b) - \sum_{ij \in E} \theta_{ij} a_{ij} - \pi \sum_{ij
\in E_{\pi}} a_{ij} - \sum_{k \in V^1} \Theta_k b_{k},
\end{align}
which is real-analytic and strictly locally convex on the open
domain $\TE$ (see \cite{ND}), as well as convex and continuously
differentiable on the open convex set $\reals^{E_T^0} \times
\reals_{+}^{E_T^1 \cup V_1}\, \supset \,\TE$ since it is a sum of
convex continuously differentiable functions minus a linear
function. The set $E_T^1 = \{ij \in E_T \, | \, i, j \in V^1\}$
and $E_T^0 = E_T \, \setminus \, E_T^1$.

The constructibility of the pattern we are after comes from the
following result, which is proved in \cite{ND} and which plays a
central role in the proof of Theorem \ref{Thm_main}.
\begin{thm} \label{thm_variational_principle}
The unique up to isometry hyper-ideal circle pattern which
realizes the data $(\Triang, \tilde{\theta}, \Theta)$, and hence
the original data $(\cellcomplex, \theta, \Theta)$, is represented
by a unique minimum, located inside $\TE,$ of the continuously
differentiable convex functional
$$\Uglobal_{\theta, \Theta} \, : \, \reals^{E_T^0} \times
\reals_{+}^{E_T^1 \cup V_1} \, \to \, \reals$$ defined by formula
(\ref{Formula_global_functional_general}). The existence of the
minimum is guaranteed if and only if the angle data $(\theta,
\Theta)$ satisfies conditions 1 to 4 of Theorem \ref{Thm_main}.
\end{thm}


With all these tools at hand, we can move on the next section in
which an algorithmic recipe for discrete uniformization is
outlined.

\section{Algorithm for discrete uniformization via hyper-ideal circle patterns}
\label{Sec_algorithm}

Start with a closed surface $S$ of genus one or greater, together
with one of the following two geometric data on it:

\smallskip
\noindent \emph{Type 1 data}.
\begin{itemize}
    \item Either a hyperbolic or Euclidean cone metric $d$ on $S$ with cone
singularities $\text{sing}(d)$ whose cone-angles are greater than
$2\pi$.
    \item A finite set of points $V$ on $S$ such that
$\text{sing}(d)=V_1 \subseteq V$ and $V_0 = V \setminus V_1$.
\end{itemize}

\smallskip
\noindent \emph{Type 2 data}.
\begin{itemize}

\item A finite topological branch cover $p \, : \, S \, \to \,
\hat{\CC}$ with ramification points $V_1 \subset S$ and branch
points $V_1(\hat{\CC}) \subset \hat{\CC}.$

\item A finite set of points $V_{\hat{\CC}} = V_0(\hat{\CC}) \cup
V_1(\hat{\CC})$ on $\hat{\CC}$ where $V_0(\hat{\CC}) \cap
V_1(\hat{\CC}) = \varnothing$.

\item A finite set $V = p^{-1}(V_{\hat{\CC}})$ on $S$ with $V_0 =
V \setminus V_1$.



\end{itemize}

\medskip
\noindent \textbf{Step 1.} Generate the Delaunay circle pattern
either on $S, d$ with respect to $V$ if given Type 1 data, or on
$\hat{\CC}$ with respect to $V_{\hat{\CC}}$ if given Type 2 data.

\medskip
\noindent \textbf{Step 2.} In the case of Type 1 data, the
Delaunay circle pattern from Step 1 gives rise to a combinatorial
cell complex $\cellcomplex=(V, E, F)$ and an angle assignment
$\theta \, : \, E \, \to \, (0,\pi)$ of intersection angles
between adjacent Delaunay circles.

\smallskip
In the case of Type 2 data, the Delaunay circle pattern from Step
1 gives rise to a combinatorial cell complex
$\cellcomplex_{\hat{\CC}}=(V_{\hat{\CC}}, E_{\hat{\CC}},
F_{\hat{\CC}})$ and an angle assignment $\hat{\theta} \, : \,
E_{\hat{\CC}} \, \to \, (0,\pi)$ of intersection angles between
adjacent Delaunay circles.

\medskip
\noindent \textbf{Step 3.} Only in the case of Type 2 data, lift
the complex $\cellcomplex_{\hat{\CC}}$ to a cell complex
$\cellcomplex=(V, E, F)$ on $S$ via the branch covering map $p$.
Thus, $p^{-1} \big( V_{\hat{\CC}} \big) = V, \,\,\, p^{-1} \big(
E_{\hat{\CC}} \big) = E, \,\,\, p^{-1} \big( F_{\hat{\CC}} \big) =
F$. Furthermore, define the lifted angle assignment $\theta \, :
\, E \, \to \, (0,\pi)$ as $\theta_{ij} = \hat{\theta}_{p(ij)}$
for all $ij \in E$.


\medskip
\noindent \textbf{Step 4.} Subtriangulate $\cellcomplex$ and
obtain the combinatorial triangulation $\Triang = (V, E_T, F_T)$
by adding a maximal number of diagonals with non-intersecting
interiors in each non-triangular face of $\cellcomplex$ (see
Section \ref{Sec_Spaces_of_triangles_and_patterns}). Define
$E_{\pi}$ as the set of all added diagonals, also called redundant
edges in Section \ref{Sec_Spaces_of_triangles_and_patterns}. Thus,
$E_T = E \cup E_{\pi}$.

\medskip
\noindent \textbf{Step 5.} Form the set $E_T^1$ of all edges from
$E_T$ both of whose endpoints are vertexes from $V_1$. Let $E_T^0
= E_T \setminus E_T^1$.

\medskip
\noindent \textbf{Step 6.} Form the functional
\begin{align} 
\Uglobal_{\theta}(a,b) = \sum_{\Delta \in F_T} \UDelta(a,b) -
\sum_{ij \in E} \theta_{ij} a_{ij} - \pi \sum_{ij \in E_{\pi}}
a_{ij} - 2\pi\sum_{k \in V} b_{k},
\end{align}
for all $(a,b)$ from the convex set $\reals^{E_T^0} \times
\reals^{E_T^1 \cup V_1}_{+}$. The functions $\UDelta$ are defined
by formula (\ref{Formula_local_functional}) relying on the real
analytic expressions for the angles $\alpha^{\Delta}_{ij}$ and
$\beta_k^{\Delta}$ in terms of $(a,b)$, given in the second half
of Section \ref{Sec_Spaces_of_triangles_and_patterns}. Recall that
the angle functions can be continuously extended as linear
functions outside the domains of their initial (real analytic)
definition (see Section \ref{Sec_Variational principle}).
Furthermore, also in Section \ref{Sec_Variational principle}, it
was commented that there exist analytic formulas for the volume
function $\Vol$ in terms of dihedral angles $(\alpha^{\Delta},
\beta^{\Delta})$.


\medskip
\noindent \textbf{Step 7.} As explained in Section
\ref{Sec_Variational principle}, the functional
$\Uglobal_{\theta}$ is convex and continuously differentiable on
$\reals^{E_T^0} \times \reals^{E_T^1 \cup V_1}_{+}$, and locally
strictly convex and real-analytic on its open subdomain $\TE$.
Find the unique minimum $(a^{\star},b^{\star}) \in \TE$ of
$\Uglobal_{\theta}$ whose existence is guaranteed by Theorems
\ref{thm_uniformization_negative_curvature} and
\ref{thm_general_uniformization}.

\medskip
\noindent \textbf{Step 8.} Compute the edge-lengths and vertex
radii $(l^{\star},r^{\star}) = \Psi(a^{\star},b^{\star})$ using
formulas (\ref{Eqn_Hyp_r_b}) to (\ref{Eqn_Hyp_l_a}).

\medskip
\noindent \textbf{Step 9.} Following the combinatorics of
$\Triang,$ lay out in $\hyperbolicplane$ the hyperbolic triangles
determined by the edge-lengths $l^* : E_T \to \reals_{+}$. Thus, a
fundamental domain of a Fuchsian group \cite{Bus,ThuNotes,ThuBook}
is obtained and if one computes the generators of this Fuchsian
group, a discrete equivalent of the classical uniformization
theorem is obtained. If desired, erase the redundant edges
$E_{\pi}$ to represent accurately the geodesic realization of the
complex $\cellcomplex$ in $\hyperbolicplane$.

\medskip
\noindent \textbf{Step 10} (\emph{Optional})\textbf{.} One could
also draw the resulting hyper-ideal circle pattern, by first
drawing all vertex circles given by $r^* : V_1 \to \reals_{+}$.
Then, the presence of the vertex circles uniquely determines the
orthogonal face-circles.

In Section \ref{Sec_computer_realizations} we present some results
from the computer implementation of the algorithm.

\section{Realization of hyper-ideal circle patterns on the Riemann sphere}
\label{Sec_realization_on_sphere}

\subsection{Constructions} We start with a topological gluing
construction. Let $\cellcomplex = (V, E, F)$ be a strongly regular
complex on the two-sphere $S^2$. Fix $k_{\infty} \in V$ and define
$\bar{D}(k_{\infty}) = \cup F_{k_{\infty}}$ to be the union of all
closed faces of $\cellcomplex$ attached to $k_{\infty}$. By strong
regularity, $\bar{D}(k_{\infty})$ is a closed topological disk
embedded in $S^2$. Let $\partial \bar{D}(k_{\infty}) = \sigma$ be
its boundary, which is composed of edges of $\cellcomplex$, and
let ${D}(k_{\infty})$ be its open interior. Next, take two copies
of $(S^2, \cellcomplex)$ and remove $D(k_{\infty})$ from both of
them. Then glue the two copies together along the two copies of
boundary $\sigma$, identifying pairs of twin edges. We obtain the
connected sum $S^2 {\#}_{\sigma} S^2 \, \cong \, S^2$ together
with a strongly regular complex $\cellcomplex_{\sigma} =
(V_{\sigma}, E_{\sigma}, F_{\sigma})$ on it. Notice that
$\cellcomplex_{\sigma}$ has a topological symmetry, which is an
involution fixing point-wise the simple closed loop $\sigma$.

Now, assume that our strongly regular complex $\cellcomplex$ on
$S^2$ comes equipped with an angle assignment $\theta \, : \, E \,
\to \, (0,\pi)$ which satisfies the conditions of Bao and
Bonahon's Theorem \ref{thm_Bao}. Then there exists a unique, up to
$\mathbb{P}SL(2,\CC)$ automorphism, hyper-ideal circle pattern on
$\hat{\CC}$ which realizes the combinatorial angle data
$(\cellcomplex, \theta)$. Let us assume that the pattern has at
least one true vertex circle, i.e. $V = V_0 \sqcup V_1$ with $V_1
\neq \varnothing$. As already explained in the Discussion on
hyper-ideal circle patterns on $\hat{\CC}$ from Section
\ref{Sec_def_and_notations}, we can choose a point $k_{\infty}$
from the interior of the vertex circle and stereographically
project the pattern on $\CC$ so that the vertex circle becomes the
ideal boundary of the hyperbolic plane and the rest of the pattern
becomes a hyper-ideal circle pattern on a convex geodesic polygon
$P$ in $\hyperbolicplane$ with combinatorics $\cellcomplex
\setminus D(k_{\infty})$. The interior angle at a vertex $i$ of
$\partial P$ is equal to $\theta_{ik_{\infty}}$ and the angle
between a geodesic edge $ij$ of $\partial P$ and the face circle
of the decorated polygon from the pattern attached to $ij$ is
$\theta_{ij}$. Take two copies of $P$ and glue them together along
$\partial P$ identifying isometrically the pairs of twin edges.
The result is a sphere $S^2_P$ together with a hyperbolic metric
with cone singularities at the vertices that were once boundary
vertices of $P$. Moreover, there is a (generalized) hyper-ideal
circle pattern on $S^2_P$ which is symmetric with respect to an
isometric involution which fixes point-wise the former boundary
$\partial P$. Inside each copy of $P$, the intersection angles
between adjacent face circles are equal to the intersection angles
from the original circle pattern on $\hat{\CC}$. Hence, for a
non-boundary edge $ij$ of a decorated polygon in one of the two
copies of $P$ the intersection angle is $\tilde{\theta}_{ij} =
\theta_{ij}$. The angle between adjacent face circles at a former
boundary edge $ij \subset
\partial P$ is $\tilde{\theta}_{ij} = 2 \theta_{ij}$. For a vertex
$i$ inside a copy of $P$ the cone angle is $\Theta_i = 2\pi$ while
for a vertex at the former boundary $\partial P$ the cone angle is
$\Theta_i = 2 \theta_{ik_{\infty}}$. Observe that by construction,
the circle pattern on $S^2_P$ has combinatorics
$\cellcomplex_{\sigma}$. Therefore, it realizes the combinatorial
angle data $(\cellcomplex_{\sigma}, \tilde{\theta}, \Theta)$.
Furthermore, it is unique, up to isometry, due to its isometric
involutive symmetry and the uniqueness of its two components $P$
guaranteed by Bao and Bonahon's Theorem. Later we will see that
uniqueness also follows from a variational principle. With these
constructions and notations at hand, we are ready to proceed to
the algorithm.

\subsection{Algorithm for realization of hyper-ideal circle
patterns on $\hat{\CC}$}

Start with a topological sphere $S^2$ together with the following
data on it:

\smallskip

\begin{itemize}

\item A finite topological branch cover $p \, : \, S^2 \, \to \,
\hat{\CC}$ with ramification points $V_1 \subset S^2$ and branch
points $V_1(\hat{\CC}) \subset \hat{\CC}.$

\item A finite set of points $V_{\hat{\CC}} = V_0(\hat{\CC}) \cup
V_1(\hat{\CC})$ on $\hat{\CC}$ where $V_0(\hat{\CC}) \cap
V_1(\hat{\CC}) = \varnothing$.

\item A finite set $V = p^{-1}(V_{\hat{\CC}})$ on $S^2$ with $V_0
= V \setminus V_1$.



\end{itemize}

\medskip
\noindent \textbf{Step 1.} Generate the Delaunay circle pattern on
$\hat{\CC}$ with respect to $V_{\hat{\CC}}$.

\medskip
\noindent \textbf{Step 2.} The Delaunay circle pattern from Step 1
gives rise to a combinatorial cell complex
$\cellcomplex_{\hat{\CC}}=(V_{\hat{\CC}}, E_{\hat{\CC}},
F_{\hat{\CC}})$ and an angle assignment $\hat{\theta} \, : \,
E_{\hat{\CC}} \, \to \, (0,\pi)$ of intersection angles between
pairs of adjacent Delaunay circles.

\medskip
\noindent \textbf{Step 3.} Lift the complex
$\cellcomplex_{\hat{\CC}}$ to a cell complex $\cellcomplex=(V, E,
F)$ on $S^2$ via the branch covering map $p$. Thus, $p^{-1} \big(
V_{\hat{\CC}} \big) = V, \,\,\, p^{-1} \big( E_{\hat{\CC}} \big) =
E, \,\,\, p^{-1} \big( F_{\hat{\CC}} \big) = F$. Furthermore,
define the lifted angle assignment $\theta \, : \, E \, \to \,
(0,\pi)$ as $\theta_{ij} = \hat{\theta}_{p(ij)}$ for all $ij \in
E$. As a result of this, there is a strongly regular complex
$\cellcomplex$ on $S^2$ together with angle assignment $\theta : E
\to (0,\pi)$.

\medskip

\noindent {Remark:} All steps from here on are independent of what
the origin of the data $(\cellcomplex,\theta)$ is, as long as it
satisfies the conditions of Theorem \ref{thm_Bao}.

\medskip
\noindent \textbf{Step 4.} Take $k_{\infty} \in V_1$, remove the
open disc $D(k_{\infty})$ from $\cellcomplex$, as described above,
and form the symmetric (connected sum) cell complex
$\cellcomplex_{\sigma} = (V_{\sigma}, E_{\sigma}, F_{\sigma})$ on
the connected sum $S^2 \#_{\sigma} S^2 \, \cong \, S^2$ over the
boundary $\sigma = D(k_{\infty})$. The vertex set is naturally
split into $V_{\sigma} = V_{\sigma,1} \sqcup V_{\sigma,0}$
inherited from the splitting of $V$.

\medskip
\noindent \textbf{Step 5.} For an edge $ij \in E_{\sigma}$, if
$ij$ does not lie entirely on $\sigma$, then define
$\tilde{\theta}_{ij} = \theta_{ij}$. If $ij$ lies on $\sigma$,
then define $\tilde{\theta}_{ij} = 2 \theta_{ij}$. In both cases,
$ij$ is also interpreted as a former edge of $\cellcomplex
\setminus D(k_{\infty})$.

\medskip
\noindent \textbf{Step 6.} For a vertex $i \in V_{\sigma}$, if $i$
is not on $\sigma$, then define $\Theta_{i} = 2 \pi$. If $i \in
\sigma$, then define $\Theta_{i} = 2 \theta_{ik_{\infty}}$. In
both cases, $i$ is also interpreted as a former vertex of
$\cellcomplex \setminus D(k_{\infty})$.

\medskip
\noindent \textbf{Step 7.} Subtriangulate $\cellcomplex \setminus
D(k_{\infty})$ and obtain the combinatorial triangulation
$\Triang_{\sigma} = (V_{\sigma}, E_T, F_T)$ by adding a maximal
number of diagonals with non-intersecting interiors in each
non-triangular face of $\cellcomplex_{\sigma}$ (see Section
\ref{Sec_Spaces_of_triangles_and_patterns}). The triangulation can
be constructed so that it inherits the involutive symmetry of
$\cellcomplex_{\sigma}$. Indeed, one can first subtriangulate
$\cellcomplex \setminus D(k_{\infty})$ and then glue together two
identical copies along $\sigma$. Define $E_{\pi}$ as the set of
all added diagonals, also called redundant edges in Section
\ref{Sec_Spaces_of_triangles_and_patterns}. Thus, $E_T =
E_{\sigma} \sqcup E_{\pi}$.

\medskip
\noindent \textbf{Step 8.} Form the set $E_T^1$ of all edges from
$E_T$ both of whose endpoints are vertexes from $V_{\sigma, 1}$.
Let $E_T^0 = E_T \setminus E_T^1$.

\medskip
\noindent \textbf{Step 9.} Form the functional
\begin{align} 
\Uglobal_{\tilde{\theta},\Theta}(a,b) = \sum_{\Delta \in F_T}
\UDelta(a,b) - \sum_{ij \in E_{\sigma}} \tilde{\theta}_{ij} \,
a_{ij} - \pi \sum_{ij \in E_{\pi}} a_{ij} - \sum_{k \in V}
\Theta_k \, b_{k},
\end{align}
for all $(a,b)$ from the convex set $\reals^{E_T^0} \times
\reals^{E_T^1 \cup V_{\sigma, 1}}_{+}$. The functions $\UDelta$
are defined by formula (\ref{Formula_local_functional}) relying on
the real analytic expressions for the angles
$\alpha^{\Delta}_{ij}$ and $\beta_k^{\Delta}$ in terms of $(a,b)$,
given in the second half of Section
\ref{Sec_Spaces_of_triangles_and_patterns}. Recall that the angle
functions can be continuously extended as linear functions outside
the domains of their initial (real analytic) definition (see
Section \ref{Sec_Variational principle}). Furthermore, also in
Section \ref{Sec_Variational principle}, it was commented that
there exist analytic formulas for the volume function $\Vol$ in
terms of dihedral angles $(\alpha^{\Delta}, \beta^{\Delta})$.


\medskip
\noindent \textbf{Step 10.} As explained in Section
\ref{Sec_Variational principle}, the functional
$\Uglobal_{\tilde{\theta}, \Theta}$ is convex and continuously
differentiable on $\reals^{E_T^0} \times \reals^{E_T^1 \cup
V_{\sigma, 1}}_{+}$, and locally strictly convex and real-analytic
on its open subdomain $\TE$. Theorem \ref{thm_Bao} guarantees the
existence of a critical point in $\TE$. By strict convexity, the
critical point is unique (and it is a minimum). Hence, this is
another justification for the uniqueness of the hyper-ideal circle
pattern we are looking to construct. Find the unique minimum
$(a^{\star},b^{\star}) \in \TE$ of $\Uglobal_{\tilde{\theta},
\Theta}$ using convex optimization.

\medskip
\noindent \textbf{Step 11.} Compute the edge-lengths and vertex
radii $(l^{\star},r^{\star}) = \Psi(a^{\star},b^{\star})$ using
formulas (\ref{Eqn_Hyp_r_b}) to (\ref{Eqn_Hyp_l_a}).

\medskip
\noindent \textbf{Step 12.} By uniqueness and symmetry of the
data, the pattern represented by $(l^{\star},r^{\star})$ has an
isometric involution and thus it splits into two identical
hyper-ideal circle patterns with convex geodesic boundary. Each of
them has a smooth hyperbolic metric, so it is realizable in
$\hyperbolicplane$. Following the combinatorics of
$\Triang_{\sigma}$ only on one side of the simple loop $\sigma$,
lay out in $\hyperbolicplane$ the hyperbolic triangles determined
by the edge-lengths $l^{\star} : E_T \to \reals_{+}$. Thus, a
hyper-ideal circle patterns with convex geodesic boundary is
obtained. If desired, erase all redundant edges from $E_{\pi}$ to
represent accurately the geodesic realization of the complex
$\cellcomplex \setminus D(k_{\infty})$ in $\hyperbolicplane$.

\medskip
\noindent \textbf{Step 14.} Using the assignment of vertex radii
$r^{\star} : V_{\sigma} \to \reals_{+}$, first draw the vertex
circles centered at the vertices of the realized geodesic complex
in $\hyperbolicplane$. Then, the presence of the vertex circles
uniquely determines the orthogonal face-circles. By adding the
circle at infinity $\partial \hyperbolicplane$ to the collection
of vertex circles, and by adding the circles representing the
hyperbolic geodesics of the boundary of the complex to the
collection of face circles, we obtain the hyper-ideal circle
pattern on $\hat{\CC}$ which realizes the combinatorial angle data
$(\cellcomplex, \theta)$.

\paragraph{Remark:} The same procedure can be carried out by
choosing a vertex $k_{\infty} \in V_0$ instead of $V_1$. Since in
this case the vertex circle at $k_{\infty}$ is a point, the
cone-metric on the doubled sphere will be Euclidean and uniquely
defined up to scaling.


\section{Numerical examples} \label{Sec_computer_realizations}

We use the numerical optimization package TAO \cite{tao-user-ref}
to perform the minimization of the functional for certain
examples. We use the BLMVM method which is a quasi-Newton method
that approximates the Hessian matrix using a fixed number of
explicit gradient evaluations. We configure the bounded
domain--see Steps 5 and 6 of Section \ref{Sec_algorithm}--using
the API for bounded minimization in the TAO application.

\begin{figure}[ht]
\centering
\raisebox{0.9cm}{\includegraphics[width=0.4\textwidth]{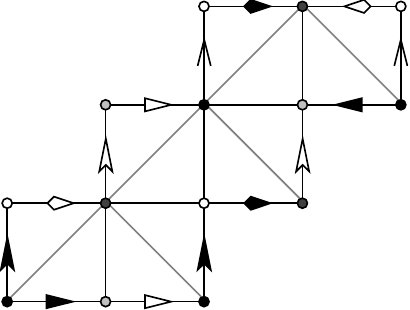}}
\hspace{0.3cm}
\includegraphics[width=0.5\textwidth]{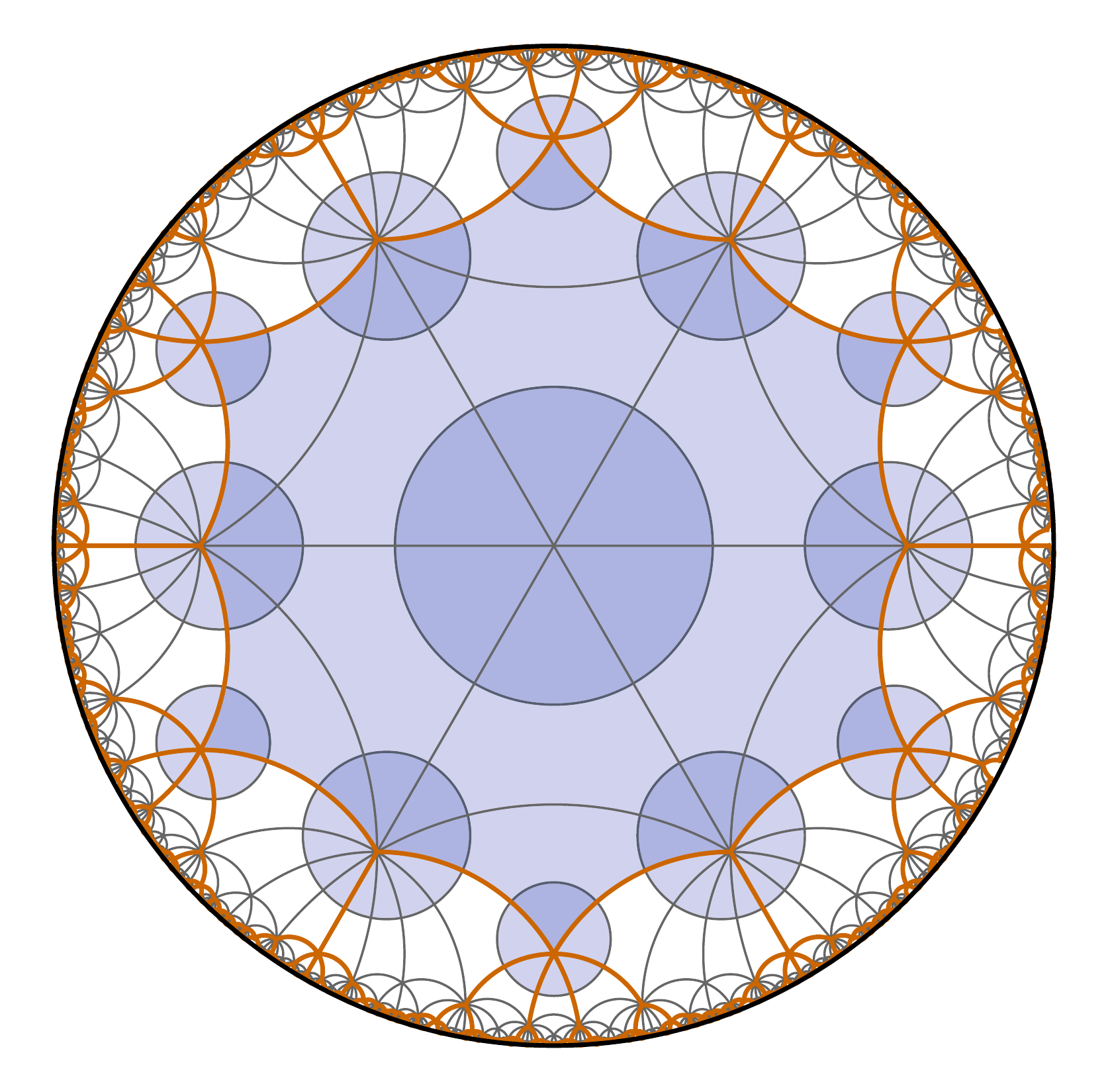}
\caption{Square-tiled representation of the Riemann surface of
Lawson's minimal surface (left). Universal cover, Fuchsian
uniformization, fundamental domain and vertex circles (right), see
Subsection~\ref{sec:examples_1_and_2}.}
\label{fig:lawson-diagonals}
\end{figure}

\subsection{Square-tiled Reimann surfaces}
\label{sec:examples_1_and_2}

\paragraph{Example 1.} As a first example we use a square-tiled surface which is
conformally equivalent to the Riemann surfaces associated to the
Lawson's genus 2 minimal surface in $\mathbb{S}^3$ \cite{Law1970},
see Figure~\ref{fig:lawson-diagonals}, left. It has four vertices
(white, light-grey, dark-grey, black), 18 edges, and 6 faces. The
6 diagonal edges are redundant edges in the sense that the face
circle intersection angle is $\pi$ at these edges, i.e. the two
circles on both sides of a redundant edge actually coincide. As a
result of this, in the hyperbolic realization Figure
\ref{fig:lawson-diagonals}, right, all pairs of decorated
triangles sharing redundant edges are merged to form decorated
quadrilaterals. Hence the conformal structure on the grey
diagonals is given by $\theta_{\mathrm{grey}}=\pi$ and on the
black edges by $\theta_{\mathrm{black}}=\frac{\pi}{2}$. Boundary
edges are identified as indicated by the arrows. In this case, all
vertices have cone angles greater than $2\pi$ so $V=V_1$ while
$V_0=\varnothing$. All of them become vertex circles after
uniformization. On Figure~\ref{fig:lawson-diagonals}, right, the
vertex circles are depicted, while the face circles are not for
the sake of better clarity of the picture.

In this example all variables $a_i$, $b_i$ are strictly positive.
We minimize the functional using the options described above. The
solver converges after 17 iterations at a solution with a gradient
norm less than $10^{-10}$.

We choose a rotationally symmetric fundamental domain where each
square is incident to a central vertex, see
Figure~\ref{fig:lawson-diagonals}, right.


\begin{figure}[ht]
\centering
\raisebox{0.9cm}{\includegraphics[draft=false,width=0.4\textwidth]{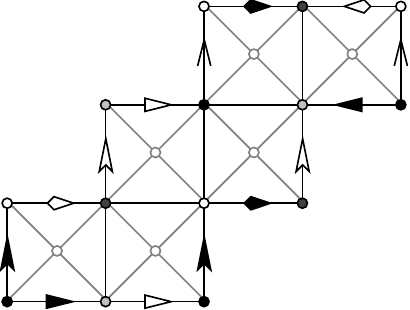}}
\hspace{0.3cm}
\includegraphics[draft=false,width=0.5\textwidth]{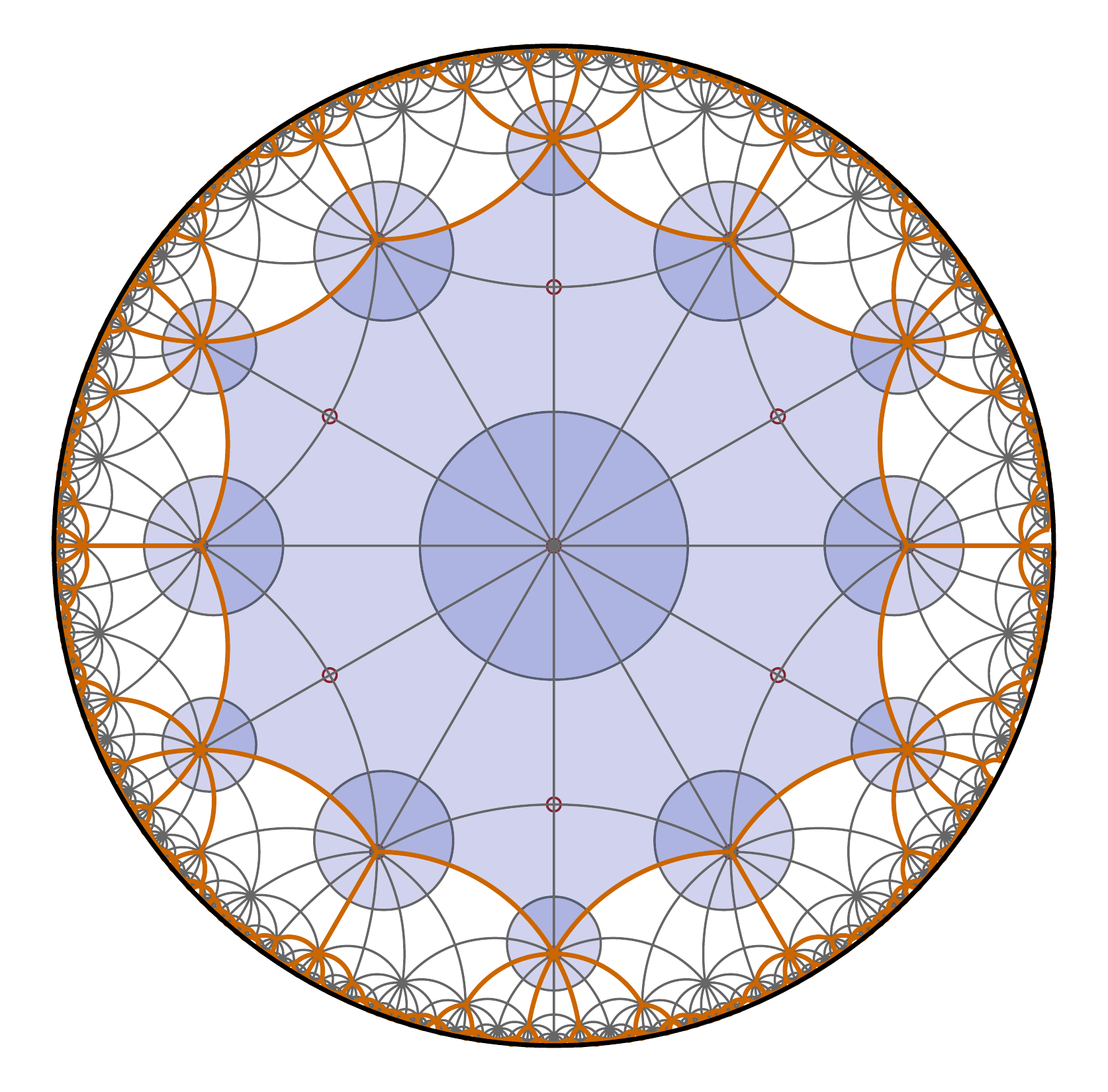}
\caption{Square-tiled representation of the Riemann surface of
Lawson's minimal surface including ideal vertices in the centers
of the squares (left). Universal cover, Fuchsian uniformization,
fundamental domain and vertex circles (right). See
Subsection~\ref{sec:examples_1_and_2}.} \label{fig:lawson_branch}
\end{figure}

\paragraph{Example 2.} For the second example, we add six new vertices to the
square-tiled surface form Figure \ref{fig:lawson-diagonals}, left.
Namely, in the center of each square we insert a vertex and
triangulate the quadrilaterals as shown on
Figure~\ref{fig:lawson_branch}, left. Note that these vertices
correspond to the ramification points of the Riemann surface when
represented as an algebraic curve in $\mathbb{C}^2$, see
Subsection~\ref{sec:examples_3_and_4}.

The conformal structure is given by $\theta_{\mathrm{black}}=\pi$
on black edges and $\theta_{\mathrm{grey}}=\frac{\pi}{2}$ on grey
edges. The black vertices are redundant, so we obtain decorated
quadrilaterals after uniformization, Figure
\ref{fig:lawson_branch}, right. As the new white vertices are flat
to begin with we exclude the corresponding variables from the
functional, i.e., $b_{\mathrm{white}} \equiv 0$. In the
corresponding hyper-ideal tetrahedron this corresponds to an ideal
vertex of the tetrahedron. At the same time when introducing ideal
vertices we change the domain of optimization for certain
variables, i.e., for edges incident with at least one ideal vertex
we have $a_{\it ij}\in \mathbb{R}$. In this case, the white
vertices form the set vertex set $V_0$ while the rest of the
vertices form the set $V_1$. As one can see on
Figure~\ref{fig:lawson-diagonals}, right, all vertices from $V_1$
become vertex circles, while the ones from $V_0$ stay points
(circles of radius zero).

Just like in the preceding example, the solver converges after 26
iterations
to a solution with a gradient norm less than $10^{-7}$. 

We choose the same fundamental domain as in the previous example
and calculate the generators of the Fuchsian uniformization group
and the corresponding universal cover, see
Figure~\ref{fig:lawson_branch}, right.

\subsection{Hyperelliptic Riemann surfaces}
\label{sec:examples_3_and_4}

\begin{figure}[ht]
\centering
\raisebox{0.8cm}{\includegraphics[draft=false,width=0.42\textwidth]{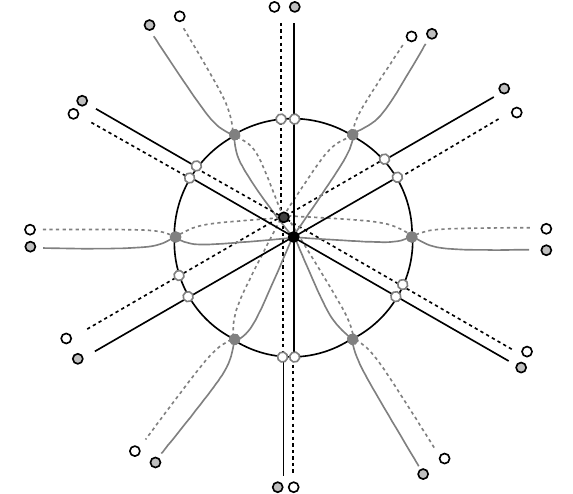}}
\hspace{0.3cm}
\includegraphics[draft=false,width=0.5\textwidth]{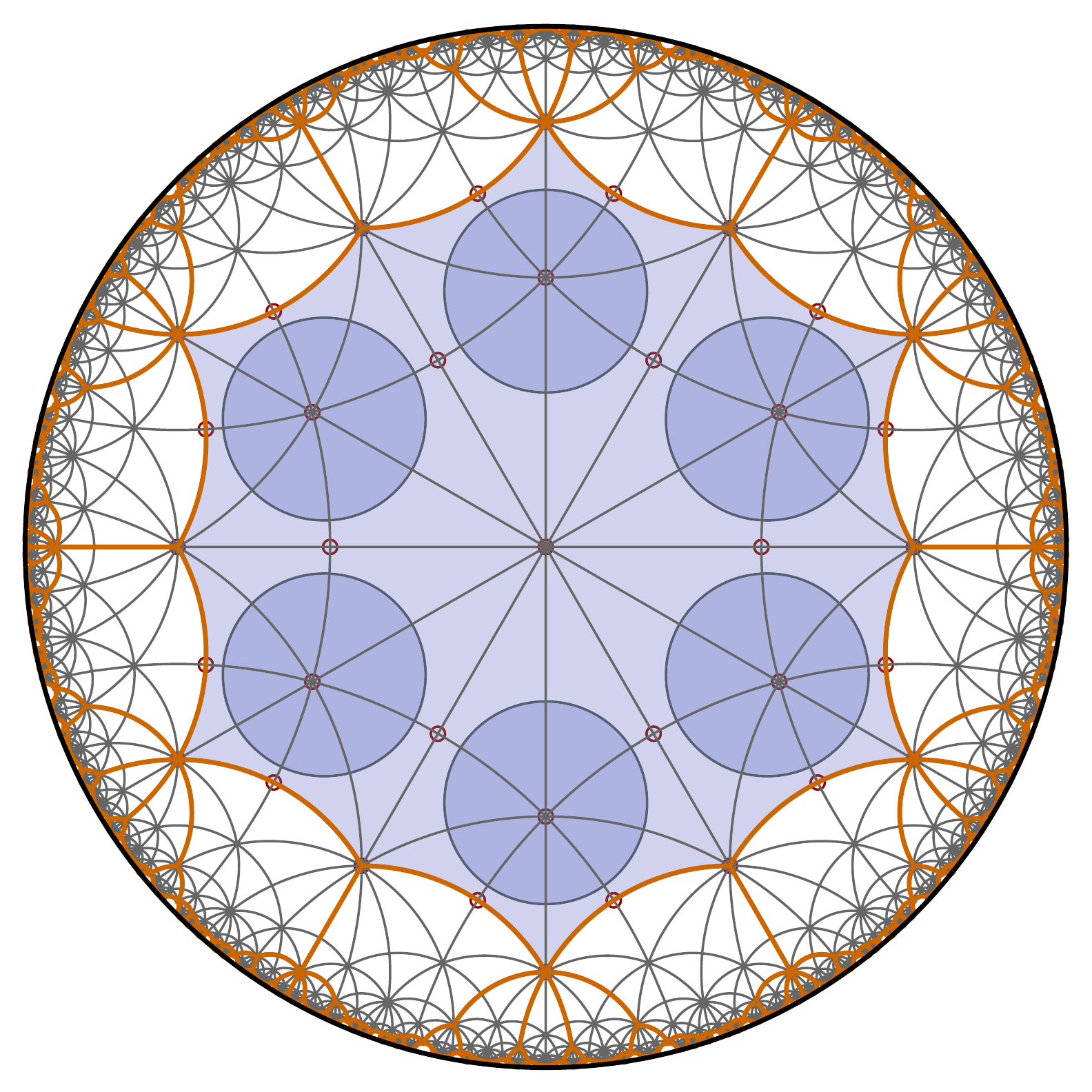}
\caption{The Riemann surface of Lawson's genus $2$ minimal surface
in $\mathbb{S}^3$ represented as doubly covered polyhedral surface
over $\hat{\mathbb{C}}$ (left). Universal cover, Fuchsian
uniformization, fundamental domain and vertex circles (right). See
Section~\ref{sec:examples_3_and_4}.} \label{fig:lawson_curve}
\end{figure}

\paragraph{Example 3.} In this example we discretely uniformize the complex algebraic
curve $\mu^2=\lambda ^6-1$. The latter is the Riemann surface
associated to Lawson's minimal surface in the three sphere
\cite{Law1970} and it is represented as a branched cover over
$\hat{\mathbb{C}}$. We generate a Delaunay triangulation on
$\hat{\CC}$ that includes the ramification points of the algebraic
curve. Hence the triangulation includes the six roots of unity as
vertices. We add the mid-points between these vertices as well as
the north and the south poles. Then we lift the corresponding
Delaunay triangulation to the algebraic curve creating a
two-sheeted cover of $\hat{\mathbb{C}}$ branched around the six
vertices at the roots of unity, see Figure~\ref{fig:lawson_curve},
left.

The conformal structure is calculated by measuring the
intersection angles of the face circumcircles in
$\hat{\mathbb{C}}\cong\mathbb{S}^2$. We first construct the
Delaunay triangulation on $\mathbb{S}^2$ using a convex hull
algorithm. Then we use a suitable stereographic projection to
measure circle intersection angles in the plane.

Using this procedure we end up with six positive variable
vertices, i.e., the vertices at branch points of the curve. The
edges are all adjacent to at least one ideal vertex. Hence all
edge variables are real, $a_{\it ij} \in \mathbb{R}$.

The solver converges after 16 iterations to an accuracy less than
$10^{-8}$ gradient norm.

The four vertices on the north and south pole correspond to the
four vertices of the quadrilaterals on the Riemann surfaces used
in the examples from Subsection \ref{sec:examples_1_and_2}. We
choose a fundamental domain with the same cuts on the surface to
produce the universal cover presented in
Figure~\ref{fig:lawson_curve}, right.


\begin{figure}[ht]
\centering
\raisebox{0.5cm}{\includegraphics[draft=false,width=0.4\textwidth]{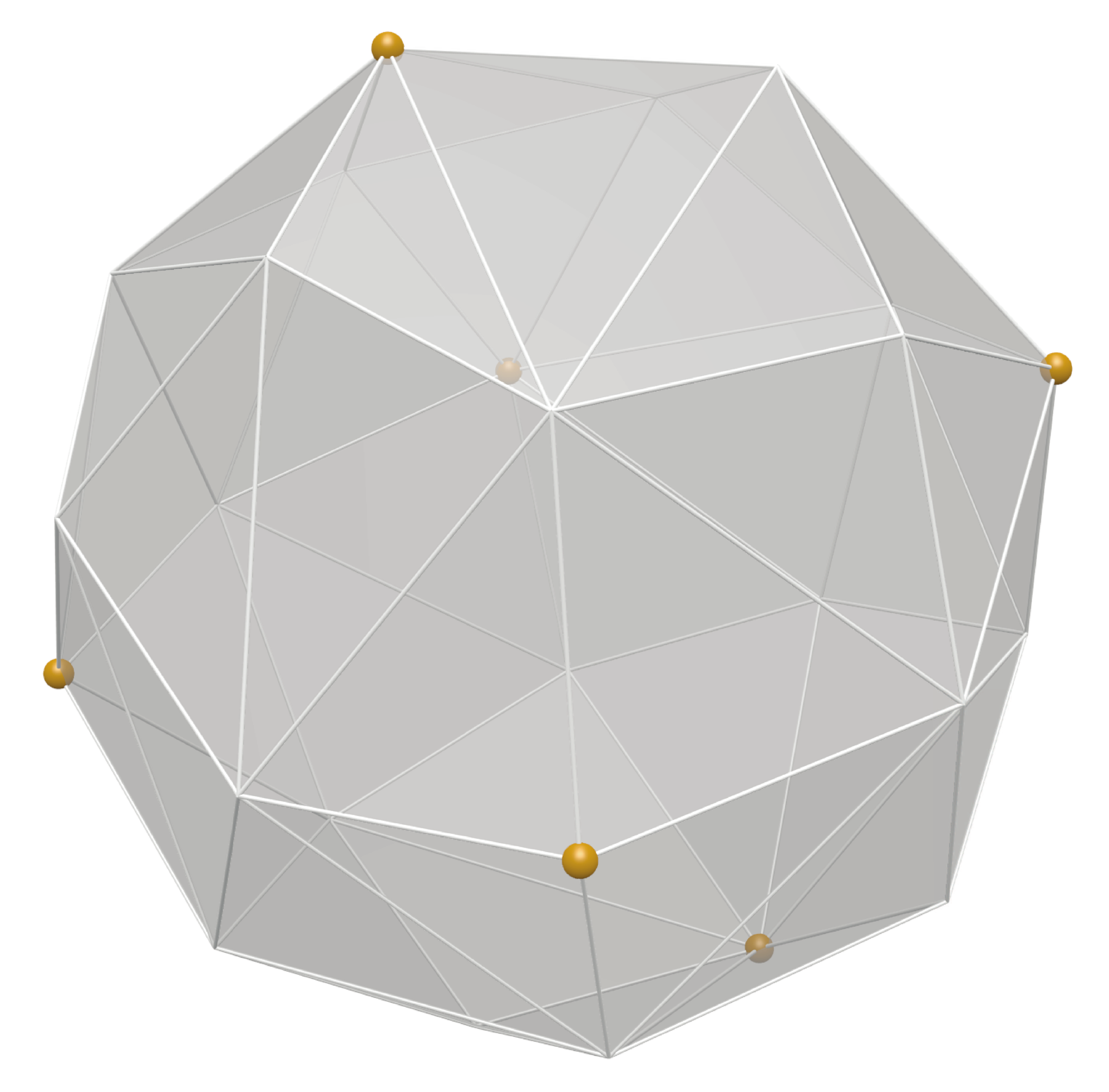}}
\hspace{0.2cm}
\includegraphics[draft=false,width=0.5\textwidth]{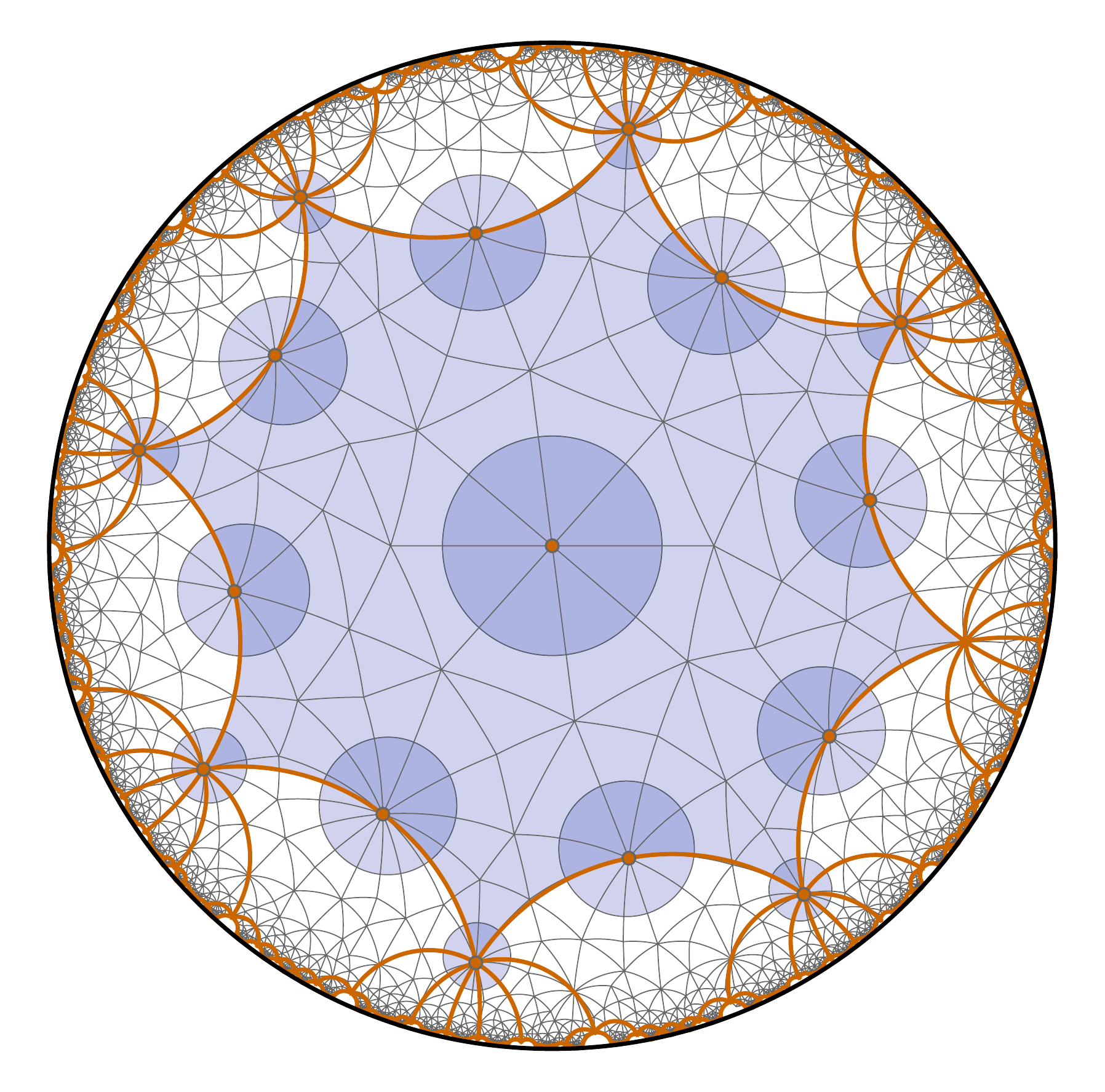}
\caption{General hyperelliptic surface of genus 2. Two-sheeted
cover of a polyhedron with vertices on the sphere (left). Large
vertices are the branch points of the corresponding algebraic
curve. Universal cover, Fuchsian uniformization, fundamental
domain and vertex circles (right). See
Section~\ref{sec:examples_3_and_4}.} \label{fig:general_curve}
\end{figure}

\paragraph{Example 4.} In this example we calculate the discrete
uniformization of a more general hyperelliptic curve of genus 2.
To achieve a better visual representation, we choose the branch
data so that the surface admits an approximately regular
fundamental domain, i.e., the branch points form approximately a
regular octahedron, see Figure~\ref{fig:general_curve}.

The surface is constructed in a way similar to the first example
of this Subsection. The triangulation includes the branch vertices
$V_1(\hat{\CC})$ and additional points $V_0(\hat{\CC})$ chosen
randomly on the sphere. Just like in the previous example, no edge
connects two branch points, hence all edge variables are real.

The conformal structure is calculated by stereographic projection.
The solver converges after 16 iteration with an accuracy of less
than $10^{-8}$.

We choose a fundamental domain that is almost a regular polygon
where opposite sides are identified. 

\subsection{Discussion}

We expect to increase the accuracy and speed of the solver if we
implement the Hessian matrix explicitly and use Newton's method as
implemented in TAO. Furthermore, by using a doubling (pillow)
construction, one can apply the so far developed methods and
algorithms to the case of a topological sphere $S$, in the spirit
of Theorem \ref{thm_sphere_uniformization}. In other words, given
the combinatorial angle data, we can construct hyper-ideal circle
patterns corresponding to the hyper-ideal polyhedra form Bao and
Bonahon's Theorem \ref{thm_Bao}.


\section*{Acknowledgements} This research is supported by DFG
(Deutsche Forschungsgemeinschaft) in the frame of
Sonderforschungsbereich/Transregio 109 ``Discretization in
Geometry and Dynamics".


\end{document}